\documentclass[10pt,twoside]{article}
\usepackage{amsfonts}
\usepackage{stmaryrd}
\usepackage{amssymb}
\usepackage{euscript}
\usepackage{pstricks}
\usepackage{amsthm}
\usepackage{latexsym,amsmath,amscd}
\usepackage{mathrsfs}
\usepackage{graphicx}
\usepackage{color}
\usepackage{dsfont}
\usepackage{enumerate}
\usepackage[normalem]{ulem}

\newtheorem{theorem}{Theorem}[section]

\newtheorem{lem}{Lemma}[section]
\newtheorem{pro}{Proposition}[section]
\newtheorem{cor}{Corollary}[section]
\newtheorem{rem}{Remark}[section]
\newtheorem{rems}{Remarks}[section]
\newtheorem{ex}{Example}[section]
\newtheorem{defi}{Definition}[section]
\newtheorem{hyp}{Assumption}[section]
\newtheorem{con}{Conjecture}[section]

\newcommand{\ssc}{\subsection}
\newcommand{\sssc}{\subsubsection}

\newcommand{\bt}{\begin{theorem}}
\newcommand{\et}{\end{theorem}}
\newcommand{\bl}{\begin{lem}}
\newcommand{\el}{\end{lem}}
\newcommand{\bp}{\begin{pro}}
\newcommand{\ep}{\end{pro}}
\newcommand{\bcor}{\begin{cor}}
\newcommand{\ecor}{\end{cor}}
\newcommand{\bcon }{\begin{con} \rm }
\newcommand{\econ }{\end{con}}
\newcommand{\lab }{\label }
\newcommand{\bd}{\begin{defi} \rm }
\newcommand{\ed}{\end{defi}}
\newcommand{\brem }{\begin{rem} \rm }
\newcommand{\erem }{\end{rem}}
\newcommand{\brems }{\begin{rems} \rm }
\newcommand{\erems }{\end{rems}}
\newcommand{\bhyp }{\begin{hyp} \rm }
\newcommand{\ehyp }{\end{hyp}}
\newcommand{\bex}{\begin{ex} \rm }
\newcommand{\eex}{\end{ex}}
\newcommand{\be}{\begin{equation}}
\newcommand{\ee}{\end{equation}}
\newcommand{\bde}{\begin{displaymath}}
\newcommand{\ede}{\end{displaymath}}
\newcommand{\beq}{\begin{eqnarray*}}
\newcommand{\eeq}{\end{eqnarray*}}
\newcommand{\ba}{\begin{align*}}
\newcommand{\ca}{\end{align*}}
\newcommand{\beqa}{\begin{eqnarray}}
\newcommand{\eeqa}{\end{eqnarray}}
\newcommand{\bel }{\left\{\begin{array}{ll}}
\newcommand{\eel}{\cr \end{array} \right.}

\newcommand{\expect}[2]{\mathbb{E}_{#1}\left( #2 \right) }
\newcommand{\cE}[3]{\mathbb{E}_{#1}\left(\left. #2\, \right| #3 \right)}
\newcommand{\cP}[3]{#1\left(\left. #2\, \right| #3 \right)}
\newcommand{\seq}[1]{{\lbrace #1 \rbrace}}

\newcommand{\llb}{\llbracket\,}
\newcommand{\rrb}{\,\rrbracket}
\newcommand{\wt}[1]{\widetilde{#1}}

\newcommand{\sbracket}[2]{[ #1, #2 ]}
\newcommand{\pbracket}[2]{\left< #1, #2 \right>}

\newcommand{\Keywords}[1]{\par\noindent{\small{\bf Keywords\/}: #1}}
\newcommand{\Class}[1]{\par\noindent{\small{\bf Mathematics Subjects Classification (2010)\/}: #1}}

\def\II{I}
\def\tII{\widetilde I}
\def\JJ{J}
\def\tJJ{\widetilde J}
\def\IIp{{}^p\II}
\def\JJp{{}^p\JJ}

\def\VV{V}

\def\BB{B}
\def\AA{A}
\def\mm{m}
\def\KK{K}
\def\LL{L}
\def\B{{\mathcal B}}

\def\IPP{{\cal I}(\P)}
\def\IPPS{{\cal I}(\P')}
\def\MM{U}
\def\whMM{U^*}

\def\UV{U}

\def\Fp{{}^pF}

\def\mp{{}^pm}

\def\wt{\widetilde}
\def\wh{\widehat }
\def\F{{\cal F}}
\def\G{{\cal G}}

\def\K{{\cal K}}
\def\ff{{\mathbb F}}
\def\gg{{\mathbb G}}
\def\hh{{\mathbb H}}
\def\kk{{\mathbb K}}
\def\rr{{\mathbb R}}
\def\err{{\bar{\mathbb R}}}
\def\nn{{\mathbb N}}
\def\Q{{\mathbb Q}}
\def\QT{\wt {\mathbb Q}}
\def\P{{\mathbb P}}

\def\EP{\mathbb{E}_{\mathbb{P}}}
\def\EPH{\mathbb{E}_{\widehat {\mathbb{P}}}}
\def\EQ{\mathbb{E}_{\mathbb{Q}}}
\def\EQT{\mathbb{E}_{\wt {\mathbb Q}}}

\def\I{\mathds{1}}

\title{{\large \bf PROGRESSIVE ENLARGEMENTS OF FILTRATIONS WITH \vskip 5 pt PSEUDO-HONEST TIMES AND THEIR APPLICATIONS \\ IN FINANCIAL MATHEMATICS}\vskip 90 pt}

\author{Libo Li and Marek Rutkowski
\footnote{The research of Marek Rutkowski was supported under Australian Research Council's Discovery Projects funding scheme
(DP120100895). The authors are grateful to Monique Jeanblanc and Shiqi Song for valuable discussions
and great hospitality during their visits to University of Evry Val d'Essonne.}
\\ School of Mathematics and Statistics
\\ University of Sydney
\\ NSW 2006, Australia}

\date{\vskip 35 pt Working paper: December 15, 2011 \vskip 15 pt}

\begin{document}
\maketitle
\vskip 45 pt
\begin{abstract}
We deal with various alternative decompositions of $\ff$-martingales with respect to the  filtration $\gg$ which represents the enlargement
of a filtration $\ff$ by a  progressive flow of observations of a random time that either belongs to the class of pseudo-honest times or satisfying
the extended density hypothesis. Several related results from the existing literature are essentially extended. Results on $\gg$-semimartingale decompositions
of $\ff$-local martingales are crucial for applications in financial mathematics, most notably in the context of modeling credit risk
and the study of insider trading where the enlargements of filtration play a vital role. We outline two potential applications
of our results to specific problems arising in financial mathematics.
\vskip 20 pt
\vskip 20 pt
\Class{60H99,$\,$91H99}

\end{abstract}


\vfill \eject

\section{Introduction} \label{sec1}

We continue here the research from Li and Rutkowski \cite{LR1}, by addressing the issues related to properties of enlarged
filtrations for various classes of random times.
We work throughout on a probability space $(\Omega, \F , \P)$ endowed with a filtration $\ff$ satisfying the usual conditions.
It is assumed throughout that $\tau $ is any $\rr_+$-valued random time given on this space.

\bd \lab{defenla}
An {\it enlargement of $\ff$ associated with} $\tau$
is any filtration $\kk = (\K_t)_{t \geq 0}$ in $(\Omega, \F , \P)$, satisfying the usual conditions,
 and such that: (i) the inclusion $\ff \subset \kk$ holds, meaning that $\F_t \subset \K_t$ for all $t \in \rr_+$,
 and (ii) $\tau $ is a $\kk$-stopping time.
\ed

Let us recall two particular enlargements, which were extensively studied in the literature (see, e.g.,
Dellacherie and Meyer \cite{DM}, Jacod \cite{JJ}, Jeanblanc and Le Cam \cite{JL}, Jeulin \cite{J1,J2,J3},
Jeulin and Yor \cite{JY1,JY2,JY3}, and Yor \cite{Y1,Y2}).
Section 8 in Nikeghbali \cite{N} provides an overview of the most pertinent results, but mostly under either
the postulate $({\bf C})$ that all $\ff$-martingales are continuous and/or
the postulate $({\bf A})$ that the random time $\tau $ avoids all $\ff$-stopping times.

\bd \lab{definit}
The {\it initial enlargement} of $\ff$ is the filtration $\gg^* = (\G^*_t)_{t \geq 0}$ where the $\sigma$-field $\G^*_t$
given by the equality $\G^*_t = \cap_{s>t}\left( \sigma (\tau )\vee \F_s\right)$ for all $t \in \rr_+$.
\ed

The initial enlargement does not seem to be well suited for a general analysis of properties of a random time with
respect to a reference filtration $\ff$ since
it implies, in particular, that $\sigma (\tau ) \subset \G^*_0$, meaning that all the information about $\tau $ is already
available at time 0 (note, however, that this feature can indeed be justified when dealing with some problems
related to the so-called {\it insider trading}).
One can argue that the following notion of the progressive enlargement of $\ff$ with observations of
a random time is more suitable for
formulating and solving various problems associated with an additional information conveyed by a random time $\tau $.

\bd \lab{defprog}
The {\it progressive enlargement} of $\ff$ is the minimal enlargement, that is, the smallest
filtration $\gg = (\G_t)_{t \geq 0}$, satisfying the usual conditions, such that $\ff \subset \gg$ and $\tau$ is a $\gg$-stopping time.
More explicitly, $\G_t = \cap_{s>t} \, \G^o_{s}$ where we denote $\G^o_t =\sigma ( \tau \wedge t) \vee \F_t$
for all $t \in \rr_+$.  
\ed

Let $\hh$ be the filtration generated by the  indicator process $H_t = \I_{\{ \tau \leq t \}}$. It is apparent
that the inclusions $\ff \vee \hh \subset \gg \subset \gg^*$ are valid and, in fact, $\gg$ coincides with
the minimal enlargement of $\ff \vee \hh$ satisfying the usual conditions.
In what follows, we will mainly work with the progressive enlargement $\gg$, although in some circumstances
we will also make use of the initial enlargement $\gg^*$.

Recall that for any two filtrations $\ff \subset \kk $ on a probability space $(\Omega, \G ,\P)$,
the {\it hypothesis} $(H')$ holds for $\ff$ and $\kk$ under $\P$ whenever any $(\P,\ff)$-semimartingale is also a $(\P, \kk )$-semimartingale
(see, e.g., Dellacherie and Meyer \cite{DM}, Jeulin \cite{J2}, Jeulin and Yor \cite{JY2} or Yor \cite{Y1}).
The problem of checking whether the hypothesis $(H')$ is satisfied and  finding the canonical semimartingale decomposition of a $(\P,\ff)$-special semimartingale with respect to a progressive
enlargement $\gg$ of a filtration $\ff$ have attracted a considerable attention and were examined in several papers during the past thirty years. In particular, the
following fundamental properties are worth to be recalled: \hfill \break
(i) a $(\P,\ff)$-semimartingale may fail to be a $(\P,\gg)$-semimartingale, in general, \hfill \break
(ii) any  $(\P,\ff)$-special semimartingale stopped at $\tau $ is a $(\P,\gg)$-special semimartingale, \hfill \break
(iii) any  $(\P,\ff)$-special semimartingale is a $(\P,\gg)$-special semimartingale when $\tau $ is an {\it honest time} with respect
to a filtration $\ff$, that is, the random variable $\tau $ is an end of some $\ff$-optional set.

Furthermore, by the classic result due to Jacod \cite{JJ}, the hypothesis $(H')$ is satisfied in the case of
the initial enlargement of $\ff$ provided that $\tau$ is an {\it initial time.} Recall that
$\tau$ is called an {\it initial time} with respect to a filtration $\ff$ if there exists a measure $\eta $ on
$(\bar\rr_+, \B(\bar\rr_+))$ such that the $(\P,\ff)$-conditional distributions of $\tau$ are absolutely continuous with respect
to $\eta $, that is, $F_{du,t} \ll \eta (du)$. This property is also frequently referred to as the {\it density hypothesis}.
In the path-breaking paper by Jeulin and Yor \cite{JY1} (see also Jeulin and Yor \cite{JY2}),
the authors derived the  $(\P,\gg)$-semimartingale decomposition of the stopped process $\MM_{\tau\wedge t}$
for any random time $\tau $ and any $(\P,\ff)$-local martingale $\MM$. They also obtained the $(\P,\gg)$-semimartingale
decomposition of an arbitrary $(\P,\ff)$-local martingale $\MM$
under an additional assumption that $\tau $ is an honest time with respect to the filtration $\ff$.
The latter result was recently extended to the case of initial times by
 El Karoui et al. \cite{EJJ}, Jeanblanc and Le Cam \cite{JL}, Kchia et al. \cite{KLP} and
Nikeghbali and Yor \cite{NY}.
For obvious reasons, we are not in a position to discuss all abovementioned papers in detail here,
although some results from them will be quoted
or referred to in what follows. Let us only mention here that, by the theorem due to Stricker \cite{S},
for arbitrary two filtrations $\ff \subset \kk$,
any a $(\P, \kk )$-semimartingale which is also an $\ff$-adapted process is necessarily a $(\P, \ff )$-semimartingale.
Therefore, in order to prove that the hypothesis $(H')$ holds for a filtration $\ff $ and its progressive
enlargement $\gg$, it suffices to show that  the hypothesis $(H')$ is satisfied by $\ff$ and any filtration $\kk$
such that $\gg \subset \kk$. A typical choice of $\kk$ in this context is the initial enlargement $\gg^*$.

The hypothesis $(H')$ should be contrasted with the stronger {\it hypothesis} $(H)$ for $\ff$ and $\kk$
under $\P$, which is also frequently referred to as the {\it immersion property} between $\ff$ and $\kk$. This hypothesis, which stipulates that any
$(\P,\ff)$-local martingale is also a $(\P, \kk )$-local martingale, was first studied in the paper by Br\'emaud and Yor \cite{BY}.
In the case of the progressive enlargement $\gg$ of a filtration $\ff$ through a random time
$\tau $ defined on the underlying probability space $(\Omega, \G ,\P)$, the immersion property for $\ff$ and $\gg$
is well known to be equivalent to the hypothesis $(H)$ introduced in Definition \ref{vfr} below
(see, e.g., Elliott et al. \cite{EJJ}) and thus no confusion may arise. Note also that the
hypothesis $(H)$, unlike the hypothesis $(H')$, is not invariant under an equivalent change of
a probability measure. However, as shown in Jeulin and Yor \cite{JY2} (Proposition 2), if the
hypothesis $(H)$ is satisfied under $\P$ by $\ff$ and an arbitrary enlargement $\kk $
and a probability measure $\Q$ is equivalent to $\P$ on ${\cal F}$ then
the hypothesis $(H')$ necessarily holds for $\ff$ and $\kk$ under $\Q$.

The main hypotheses examined in the present work are the {\it hypothesis $(H\!P)$} and
the {\it extended density hypothesis} (the {\it hypothesis} $(E\!D)$, for short), as specified in Definitions \ref{vfr} and \ref{densityr}, respectively. The corresponding classes of random times are termed {\it pseudo-honest times} and {\it pseudo-initial times}. The hypothesis $(H\!P)$ is clearly weaker than the hypothesis $(H)$ and it is known to hold, in particular,
when a random time is constructed using the {\it multiplicative approach}
(see Li and Rutkowski \cite{LR1}), as well as for the alternative construction of a random time developed in Jeanblanc and Song \cite{JS1}. It was also shown in \cite{LR1} that, under mild technical assumption, the hypothesis $(H\!P)$ is equivalent to the {\it separability} of the $(\P,\ff)$-conditional distribution of $\tau $ (see Definition \ref{comsep}). The hypothesis $(E\!D)$ extends the density hypothesis; it is introduced in order to avoid the awkward assumption
on strict positivity of $(\P,\ff)$-conditional distribution of the random time $\tau$.
It is worth to point out that most results obtained for initial times can be extended to this new setting.

The paper is organized as follows. In Section \ref{sec2}, we recall some  basic properties of $(\P,\ff)$-conditional distributions of random times
and enlarged filtrations.  In particular, we provide an alternative characterization of the progressive enlargement $\gg$.
This alternative characterization of $\gg$ is used in Section \ref{sec3} in computations of conditional expectation
of $\gg$-adapted processes under the hypothesis $(H\!P)$ and the extended density hypothesis. Subsequently,
in Section \ref{sec4}, we provide sufficient conditions for a $\gg$-adapted process to be a $(\P,\gg)$-martingale.
Explicit computations of the $\gg$-compensator (that is, the $(\P,\gg)$-dual predictable projection) of the indicator process
$H_t = \I_\seq{\tau \leq t}$ are provided in Section \ref{sec5}.
Main results of this paper are established in Section \ref{sec6} in which the validity of the hypothesis $(H')$ is studied
for the progressive enlargement of the underlying filtration $\ff$ through either a pseudo-honest or a
pseudo-initial random time. We extend there several related results from the existing literature.
First, in Theorem \ref{piu}, we compute a general semimartingale decomposition of a $(\P,\ff)$-martingale
with respect to the progressively enlarged filtration $\gg$ when $\tau $ is assumed to be a pseudo-honest time.
Particular examples of this decomposition are subsequently examined in \ref{sec62} in which we postulate that a random time
was constructed using the multiplicative approach developed in \cite{LR1}, that is, using either a predictable or an optional multiplicative
system associated with a given in advance Az\'ema submartingale $F$. Finally, in Section \ref{sec63}, we deals with the corresponding results
for a pseudo-initial time. It is worth stressing that results on a $(\P,\gg)$-semimartingale decomposition of a $(\P,\ff)$-local martingale are crucial for applications
in financial mathematics, especially in credit risk models, where a random time $\tau $ represents the moment
of occurrence of some credit event (e.g., a default event).
Two examples of applications of our results to problems of financial mathematics are outlined in Section \ref{sec7}.

\vfill \eject

\section{Random Times and Filtrations} \label{sec2}

In this section, we deal with the most pertinent properties of random times and the associated enlargements of a reference filtration $\ff$.
For more details, we refer to \cite{LR1} where, in particular, various constructions of a random time are examined.
The interested reader may also consult papers by  Jeanblanc and Song \cite{JS1,JS2} for closely related results.

\subsection{Properties of Conditional Distributions} \lab{xsect4}

Let us first introduce the notation for several pertinent characteristics of a finite random time
$\tau$ defined on a filtered probability space $(\Omega , \F , \ff , \P)$.
The $(\P,\ff)$-supermartingale $G_t = \P (\tau > t\,|\, \F_t )$ is commonly
known as the {\it Az\'ema supermartingale} of $\tau $.
We will sometimes refer to the $(\P,\ff)$-submartingale $F=1-G$ as the {\it Az\'ema submartingale} of $\tau $.
The \emph{$(\P ,\ff )$-conditional distribution of} $\tau $ is the random field $(F_{u,t})_{u,t\in \err_+}$ given by
\be \lab{tt6}
F_{u,t} = \P (\tau\leq u \,|\, \F_t) , \quad \forall \, u,t \in \err_+.
\ee
The following definition characterizes the class of all conditional distributions of a random time.

\bd \label{hyp2}
A random field $(F_{u,t})_{u,t \in \err_+}$ on a filtered probability space $(\Omega , \F , \ff , \P)$ is said to be an {\it $(\P,\ff)$-conditional
distribution} if it satisfies: \hfill \break
(i) for every $u \in \err_+$ and $t \in \err_+$, we have $0\leq F_{u,t}\leq 1$, $\P$-a.s., \hfill \break
(ii) for every $u\in \err_+$, the process $(F_{u,t})_{t\in \err_+}$ is a $(\P,\ff)$-martingale, \hfill \break
(iii) for every $t\in \err_+$, the process $(F_{u,t})_{u\in \err_+}$ is right-continuous,
increasing and $F_{\infty,t} = 1$.
\ed

Note that for every $u\in \err_+$, conditions (i)-(ii) in Definition \ref{hyp2} imply that $F_{u,\infty} = \lim_{\,t \to \infty}
F_{u,t} $ and $F_{u,t} = \EP (F_{u,\infty }\,|\, \F_t)$ for every $t\in \err_+$.
Since (iii) yields $F_{u,t}\leq F_{s,t}$ for all $u \leq s$, the (non-adapted) process $(F_{u,\infty})_{u \in \err_+}$ is increasing
and thus it admits a c\`adl\`ag version. It is known (see, e.g., \cite{LR1}) that for any random field $(F_{u,t})_{u,t \in \err_+}$ there exists a random
time $\tau $ on an extension of $(\Omega , \F , \ff , \P)$ such that (\ref{tt6}) holds.

Let us examine some pertinent properties of conditional distributions of random times.
Throughout this section, by a $(\P,\ff)$-conditional distribution, we mean any random field  $(F_{u,t})_{u,t \in \err_+}$ satisfying
Definition \ref{hyp2}.  We first recall the classic hypothesis $(H)$, which was studied in numerous papers (see, e.g., Br\'emaud and Yor~\cite{BY}
or Elliott et al.~\cite{EJY}), and its generalization termed the hypothesis $(H\!P)$
(it is obvious that the hypothesis $(H)$ implies $(H\!P)$).

\bd \label{vfr}
A $(\P,\ff)$-conditional distribution $(F_{u,t})_{u,t \in \err_+}$ is said to satisfy: \hfill \break
(i) the {\it hypothesis $(H)$} whenever for all $0 \leq u \leq s < t$
\be \label{hprop}
F_{u,s} = F_{u,t},
\ee
(ii) the {\it hypothesis $(H\!P)$} whenever for all $0 \leq u < s < t$
\be \label{hpx}
F_{u,s}F_{s,t} = F_{s,s} F_{u,t}.
\ee
\ed

It was shown in \cite{LR1} that any honest time satisfies the hypothesis $(H\!P)$ and, in fact,
an $\F_\infty$-measurable random time $\tau$ is an honest time if and only if it satisfies the hypothesis $(H\!P)$.
This motivates us to say that a random time is a {\it pseudo-honest time} with respect to $\ff$ whenever the
$(\P,\ff)$-conditional distribution of $\tau $ satisfies the hypothesis $(H\!P)$.

\brem \label{drivingM}
Let us observe that if $F_{u,t}$ satisfies the hypothesis $(H\!P)$ then, for all $0\leq  u \leq s\leq t$,
\be \label{hpro}
\frac{F_{u,s}}{F_{s,s}} \, F_{s,t} = F_{u,t}.
\ee
Note that the inclusion $\seq{F_{s,s} = 0} \subset \seq{F_{u,s} = 0}$ is valid  for all $0\leq  u \leq s$ and, by convention, $0/0 = 0$.
More generally, the inclusion  $\{F_{u,s} = 0\} \subset \{F_{u,t} = 0\}$ is known to hold for all $u \leq s \leq t$ (see the proof of Lemma 2.2 in \cite{LR1}).
\erem

Let us recall the concept of {\it separability} property of a $(\P,\ff)$-conditional distribution (see \cite{JS1,LR1}).

\bd \label{defsep}
We say that a $(\P,\ff)$-conditional distribution $(F_{u,t})_{u,t \in \err_+}$ is {\it completely separable} if there exists
a positive, $\ff$-adapted, increasing process $\KK$ and a positive $(\P,\ff)$-martingale
$\LL$ such that $F_{u,t}= \KK_u \LL_t$ for every $u,t\in \rr_+$ such that $0 \leq u \leq t$.
\ed

It is easily seen that the complete separability of $F_{u,t}$  implies that the hypothesis~$(H\!P)$ holds.
Indeed, we have that $F_{u,s}F_{s,t}= (\KK_u \LL_s)(\KK_s \LL_t)  = (\KK_s \LL_s)(\KK_u \LL_t) = F_{s,s}F_{u,t}$ for all $0 \leq u < s <t$.
It appears, however, that the property of complete separability is too restrictive, since it does
not cover all cases of our interest. This motivates the weaker concept of {\it separability} (termed partial separability in \cite{JS1}).

\bd \lab{comsep}
We say that a $(\P,\ff)$-conditional distribution $(F_{u,t})_{u,t \in \err_+}$ is {\it separable at $v \geq 0$} if
there exist a positive $(\P,\ff)$-martingale  $(\LL^v_t)_{t \in \rr_+ }$ and a positive,
$\ff$-adapted, increasing process  $(\KK^v_u)_{u \in [v,\infty )}$
such that the equality $F_{u,t}= \KK^v_u \LL^v_t$ holds for every $v \leq u \leq t $. A $(\P,\ff)$-conditional
distribution $F_{u,t}$ is called {\it separable} if it is separable at all $v >0$.
\ed

\brem \lab{uiui}
It is known that if the $(\P,\ff)$-conditional distribution of $\tau$ is separable and $F_0 = 0$ then the hypothesis $(H\!P)$ holds,
and thus $\tau$ is a pseudo-honest time. Conversely, if $\tau $ is a pseudo-honest time and its $(\P,\ff)$-conditional distribution  $(F_{u,t})_{u,t \in \err_+}$ is non-degenerate then the random field $F_{u,t}$ is separable. For proofs of these properties and more details, the interested
reader is referred to~\cite{LR1}.
\erem

The next definition proposes an extension of the {\it density hypothesis}, which was introduced by Jacod~\cite{JJ}
and subsequently studied by numerous authors (see, e.g., \cite{EJJ,JL}).
Since random times satisfying the density hypothesis are called {\it initial times}, we find it natural to say that
a random time is an {\it pseudo-initial time}  when it satisfies Definition \ref{densityr}.
Let us recall from \cite{JJ} that the hypothesis $(H')$ is known to hold for the initial (and thus also the progressive) enlargement of $\ff$ with
an initial time.

\bd \label{densityr}
A $(\P,\ff)$-conditional distribution $F_{u,t}$  is said to satisfy the {\it extended density hypothesis}
(or, briefly, the {\it hypothesis} $(E\!D)$) if there exists a random field $(\mm_{s,t})_{s \geq 0,\, t \geq s}$
and an $\ff$-adapted, increasing process $D$ with $D_{0-}=0$ and such that, for all $0 \leq u \leq t$,
\be \label{qaz}
F_{u,t} =\int_{[0,u]} \mm_{s,t}\,dD_s
\ee
and, for every $s \in \rr_+$, the process  $(\mm_{s,t})_{t \geq s}$ is a positive $(\P ,\ff)$-martingale.
\ed

\brem
It is worth noting that the complete separability is a special case of the extended density hypothesis;
it is enough to take $m_{s,t}=L_t$ and $D_t = \KK_t$.
\erem

As one might guess, the results obtained under the extended density hypothesis are similar to those proven under the usual density hypothesis,
that is, for initial times. Nevertheless, it is convenient to
introduce it here, since it will allow us to circumvent an awkward non-degeneracy condition of the $(\P,\ff)$-conditional
distribution of a random time, which will be needed, for instance, in the proof of Proposition \ref{piu1}.
Moreover, it is worth noting that the extended density hypothesis is also satisfied when a pseudo-honest time
is constructed through the multiplicative construction, as shown in Remark 2.1 in \cite{LR2} (for a special case,
see also Theorem 5.2 in Jeanblanc and Song \cite{JS1}). Therefore, the study of pseudo-initial times is related to
our main goal which is to examine the properties of pseudo-honest times.

\bex
Let $X$ be a positive $\F_\infty$-measurable random variable and $\Lambda$ a continuous, $\ff$-adapted, increasing process.
Assume that $\xi$ is a random variable uniformly distributed on $[0,1]$ and independent of $\F_\infty$.
Let us define the random time $\tau$ by setting
\bde
\tau = \inf \big\{ u \geq 0 \,:\, 1-e^{-X\Lambda_u} > \xi \big\}.
\ede
Then we obtain, for all $u \leq t$,
\begin{align*}
\cP{\P}{\tau \leq u}{\F_t} & =  \cE{\P}{1-e^{-X\Lambda_u}}{\F_t}\\
			   & = \cE{\P}{\int^u_0 X e^{-X\Lambda_s} d\Lambda_s }{\F_t}\\
		           & = \int^u_0 \cE{\P}{X e^{-X\Lambda_s }}{\F_t} d\Lambda_s = \int^u_0 m_{s,t}\, dD_s
\end{align*}
where we define $m_{s,t} := \cE{\P}{X e^{- X \Lambda_s}}{\F_t}$ for all $s \leq t$ and we set $D = \Lambda$.
\eex

\brem \lab{remvv}
Assume that the random time $\tau$ is constructed through the multiplicative approach (see Lemma 5.1 in \cite{LR1}).
Then, using also Proposition 2.3.1 and Corollary 2.3.1 in \cite{L}, we deduce that the $(\P,\ff)$-conditional distribution
$F_{u,t}$  admits the following integral representation, for all $u \leq t$
\bde
F_{u,t}  = F_{t} -  \int_{(u,t]}\frac{F_tC_{s,t}}{{}^pF_{s}}\, dB_s = F_t - \int_{[\ell_t,t]}\frac{F_tC_{s,t}}{{}^pF_{s}}\, dB_s
+ \int_{[\ell_t,u]}\frac{F_t C_{s,t}}{{}^pF_{s}}\, dB_s = \int_{[\ell_t,u]}\frac{F_t C_{s,t}}{{}^pF_{s}}\, dB_s
\ede
where $ C_{s,t}$ is a multiplicative system associated with $F$,
the process ${}^pF$ is the $(\P,\ff)$-predictable projection of $F$ and the $\F_t$-measurable random time $\ell_t$ equals
\be \label{ell1}
\ell_t  = \sup\big\{ 0 \leq s \leq t :\,  C_{s,t} = 0\big\} = \sup\big\{ 0 \leq s \leq t :\,  C_{s-,t} = 0\big\}
\ee
where, by convention, $\sup \, \emptyset = 0$.
If the Az\'ema supermartingale $G$ is strictly positive then $\ell_t = 0$ for all $t\geq 0$ and thus $F_{u,t}$ satisfies
the hypothesis $(E\!D)$. This feature of the $(\P,\ff)$-conditional distribution obtained through the
multiplicative construction can be seen as an alternative motivation for Definition \ref{densityr}.
\erem

\subsection{Enlargements of Filtrations} \lab{sect4}

We will now analyze the basic properties of various enlargements of $\ff$ associated
with a random time $\tau $. When studying semimartingale decompositions of processes stopped at $\tau $,
it is common to use, at least implicitly, the following concept, formally introduced by Guo and Zeng \cite{GZ}.

\bd \lab{dr5}
 An enlargement $\kk$ of a filtration $\ff$ is said to be {\it admissible before} $\tau$ if the equality $\K_t \cap \seq{\tau >t }
 = \F_t \cap\seq{ \tau > t }$ holds for every $t \in \rr_+$.
\ed

In the case of a general (i.e., not necessarily honest) random time, we find it convenient to introduce the following notion,
stemming from a remark in Meyer \cite{M4}. Recall that the initial enlargement $\gg^*$ was
 introduced in Definition \ref{definit}.

\bd
The family $\wh \gg = (\wh \G_t)_{t\in \rr_+ }$ is defined by setting, for all $t \in \rr_+$,
\bde
\wh \G_t = \seq{ A \in \G \,|\,\exists   A_t \in \F_t \ \mbox{and} \ A^*_{t} \in \G^*_t \
\mbox{such that}\ A  = (  A_t \cap\seq{ \tau > t }) \cup ( A^*_{t}\cap \seq{\tau \leq t} ) }.
\ede
\ed

We note that, for all $t \in \rr_+$,
\be \lab{cxz}
\wh \G_t \cap\seq{ \tau > t } = \F_t \cap\seq{ \tau > t },
\quad  \wh \G_t \cap \seq{\tau \leq t} =  \G^*_t \cap \seq{\tau \leq t}.
\ee
It can be checked that the $\sigma$-field $\wh \G_t $ is uniquely characterized by conditions (\ref{cxz}).
The next elementary result shows that the family $\wh \gg$ coincides in fact with the progressive enlargement $\gg$, which was
introduced in Definition \ref{defprog} (for the proof of the lemma, we refer to \cite{LR2}).

\bl \lab{gtr}
For any random time $\tau $ the progressive enlargement $\gg$ coincides with the filtration~$\wh \gg $.
\el

\begin{proof}
Recall that $\G_t = \cap_{s>t} (\sigma ( \tau \wedge s) \vee \F_s )$  and $\G^*_t = \cap_{s>t} ( \sigma (\tau )\vee \F_s )$.
To show that $\wh \G_t = \G_t$, it suffices to check that conditions (\ref{cxz}) are  satisfied
by $\G_t$. The following relationship for all $t \in \rr_+$ is immediate
\bde
\F_t \cap\seq{ \tau > t } \subset \G_t \cap\seq{ \tau > t } \subset\wt \G_t \cap \seq{\tau > t} =  \F_t \cap \seq{\tau >t }.
\ede
This shows $\G_t \cap\seq{ \tau > t } =\F_t \cap\seq{ \tau > t }$, while also
\bde
\G_t \cap \seq{\tau \leq t} = \cap_{s>t} (\sigma ( \tau \wedge s) \vee \F_s ) \cap \seq{\tau \leq t}
= \cap_{s>t} ( \sigma (\tau )\vee \F_s )\cap \seq{\tau \leq t} = \G_t^*\cap\seq{\tau \leq t}
\ede
since $\sigma (\tau \wedge s )\cap \seq{\tau \leq t} = \sigma (\tau ) \cap \seq{\tau \leq t}$ for every $s>t$.
\end{proof}

It is easy to see that the filtration $\wh \gg$ is admissible before $\tau $.
When dealing with a semimartingale decomposition of an $\ff$-martingale after $\tau $ we will
use the following definition.

\bd
We say that an enlargement $\kk$ is {\it admissible after} $\tau$ if the equality
$\K_t \cap \{ \tau \leq t \} = \G^*_t \cap\seq{\tau\leq t}$ holds for every $t \in \rr_+$.
\ed

It is clear that the filtration $\wh \gg $ (and thus also $\gg$) is admissible after $\tau $ for any random time.
Note also that if an enlargement $\kk$ is admissible before and after $\tau $ then necessarily $\kk = \gg$.
For the proof of the next elementary lemma, we refer to \cite{LR2}.

\bl \label{gcond}
For any integrable, $\G$-measurable random variable $X$ and any enlargement $\kk = (\K_t)_{t\geq 0}$
admissible after $\tau$ we have that, for any $t\in \rr_+$,
\be \lab{t6t6}
\cE{\P}{\I_\seq{\tau\leq t} X}{\K_t} =\lim_{s\downarrow t}\, \cE{\P}{\I_\seq{\tau\leq t}X}{\sigma(\tau)\vee \F_s}.
\ee
\el

\begin{proof}
It suffices to show that
\be \lab{t6t66}
\cE{\P}{\I_\seq{\tau\leq t} X}{\K_t} = \cE{\P}{\I_\seq{\tau\leq t}X}{\G^*_t}.
\ee
The second equality in (\ref{t6t6}) will then follow from Corollary 2.4 in \cite{RY} since $\G^*_t = \cap_{s>t}\left( \sigma (\tau )\vee \F_s\right)$.
To establish (\ref{t6t66}), we will first check that, for every $A\in \K_t$,
\bde
\expect{\P}{\I_A\I_\seq{\tau\leq t} X} =\expect{\P}{\I_A\I_\seq{\tau\leq t}\cE{\P}{X}{\G^*_t}}.
\ede
Since, by assumption, $\K_t\cap\seq{\tau\leq t} = \G^*_t \cap\seq{\tau\leq t}$,
there exists an event $B \in \G^*_t$ such that $A\cap\seq{\tau \leq t} = B\cap\seq{\tau \leq t}$.
Consequently,
\bde
\EP \big( \I_A\I_\seq{\tau\leq t} X \big)  =  \EP \big( \I_B\I_\seq{\tau\leq t} X \big)
 =  \EP \big( \I_B\I_\seq{\tau\leq t} \EP (X\,|\, \G^*_t) \big)
=  \EP \big( \I_A\I_\seq{\tau\leq t}  \EP (X\,|\, \G^*_t) \big).
\ede
Hence
\bde
\EP ( \I_\seq{\tau\leq t} X\,|\, \K_t ) =  \EP \big( \I_\seq{\tau\leq t}  \EP (X\,|\, \G^*_t) \,|\,  \K_t\big) =\I_\seq{\tau\leq t} \EP (X\,|\, \G^*_t),
\ede
since the random variable $\I_\seq{\tau\leq t} \cE{\P}{X}{\G^*_t}$ is $\K_t$-measurable.
\end{proof}

\section{Conditional Expectations under Progressive Enlargements} \label{sec3}

In the rest of the paper, we work under the assumption that the $(\P,\ff)$-conditional distribution of a random time $\tau$ satisfies either
the hypothesis $(H\!P)$ or the  hypothesis $(E\!D)$, which were introduced in Definitions \ref{vfr} and \ref{densityr}, respectively.
In addition, the special case of the complete separability will be examined as well.
We will need the following auxiliary result, which
ensures that the processes  $\I_\seq{\tau >t}(G_t)^{-1}$ and $\I_\seq{\tau \leq t}(F_t)^{-1}$ are well defined
(for its proof, see \cite{LR2}).

\bl \lab{newnn}
The following inclusions hold, for every $t \in \rr_+$:
(i) $\{\tau >t\} \subset \{ G_t>0\}$, $\P$-a.s., and
(ii) $\{\tau \leq t\} \subset \{ F_t>0\}$, $\P$-a.s.
\el

\proof
Let us denote $A = \{ F_t = 1\} =  \{ \P (\tau \leq t\,|\, \F_t ) = 1\}$. Since $A \in \F_t$
\bde
 \P (A) = \int_A F_t \, d\P = \int_A \P (\tau \leq t\,|\, \F_t ) \, d\P =  \int_A \I_{\{ \tau \leq t \}} \, d\P
 = \P ( A \cap \{ \tau \leq t \}).
\ede
Hence $A= \{ F_t = 1\} = \{ G_t = 0 \}  \subset \{ \tau \leq t \}$, $\P$-a.s., and thus  $\{\tau > t\} \subset \{ G_t>0\}$, $\P$-a.s.
For part (ii), let us denote $B = \{ G_t = 1\} =  \{ \P (\tau > t\,|\, \F_t ) = 1\}$. Since $B \in \F_t$
\bde
 \P (B) = \int_B G_t \, d\P = \int_B \P (\tau > t\,|\, \F_t ) \, d\P =  \int_B \I_{\{ \tau > t \}} \, d\P
 = \P ( B \cap \{ \tau > t \}).
\ede
Hence $B= \{ G_t = 1\} = \{ F_t = 0 \}  \subset \{ \tau > t \}$, $\P$-a.s., and thus  $\{\tau \leq t\} \subset \{ F_t>0\}$, $\P$-a.s.
\endproof

\brem
Part (i) in Lemma \ref{newnn} can also be demonstrated as follows.
Let $\tau_0 = \inf \, \{ t \in \rr_+:\, G_t = 0 \ \mbox{or}\ G_{t-}=0\}$. Since $G$ is a supermartingale,
it is equal to zero after $\tau_0$ and thus
\bde
\P ( \tau_0 < \tau ) = \EP ( \I_{\{ \tau_0 < \tau \}}) = \EP ( G_{\tau_0}  \I_{\{ \tau_0 < \infty \}} ) =0.
\ede
This in turn implies that $\{\tau >t\} \subset \{ G_t>0\}$, $\P$-a.s.
\erem

\brem
Let us set, by convention, $0/0=0$. Hence, by Lemma \ref{newnn}, the quantities $\I_{\{\tau >t\}} G^{-1}_t$
and $\I_{\{\tau \leq t\}} F^{-1}_t$ are well defined for all $t$, $\P$-a.s.
\erem

\ssc{Conditional Expectations for Pseudo-Honest Times}

For a fixed $T>0$, we consider the map $\UV_T : \rr_+ \times \Omega \to \rr $ and we use the notation
$(u,\omega ) \mapsto \UV_{u,T} (\omega )$.  We postulate that $ U_T$ is a ${\cal B}(\rr_+) \otimes \F_T$-measurable map,
so that $\UV_{\tau ,T}$ is a $\sigma (\tau ) \vee \F_T$-measurable random variable. The following result corresponds to
Theorem 3.1 in El Karoui et al. \cite{EJJ}, where the case of the density hypothesis was studied.

\bl \lab{lmu1}
Let $\UV_{\cdot ,T} : \rr_+ \times \Omega \to \rr $ be a ${\cal B}(\rr_+) \otimes \F_T$-measurable map.
Assume that $\tau$ is a pseudo-honest time and the random variable $\UV_{\tau ,T}$ is $\P$-integrable. Then: \hfill \break
(i) For every $t \in [0,T)$, we have that
\bde 
\EP (  \UV_{\tau ,T} \,|\, \G_t ) =  \I_\seq{\tau >t} \wt \UV_{t,T} + \I_\seq{\tau\leq t}\wh \UV_{\tau ,t, T}
\ede
where
\be \label{rt1}
\wt \UV_{t,T} = (G_t)^{-1} \, \EP \big(   \I_\seq{\tau >t} \UV_{\tau ,T} \,|\, \F_t \big)
= (G_t)^{-1} \, \EP \bigg(  \int_{(t,\infty ]} \UV_{v,T} \, dF_{v,T} \,\Big|\, \F_t \bigg)
\ee
and, for all $0 \leq u \leq t < T$,
\be \label{xrt2}
\wh \UV_{u,t,T} = (F_t)^{-1} \, \EP (  F_{t,T} \UV_{u,T} \,|\, \F_t ).
\ee
(ii) If, in addition, $F_{u,t}$ is completely separable so that $F_{u,t}=K_uL_t$ for $u \leq t$ then (\ref{rt1}) yields
\bde 
\EP \big(   \I_\seq{T \geq \tau >t} \UV_{\tau ,T} \,|\, \F_t \big) =  \EP \bigg(  L_T \int_{(t,T]} \UV_{v,T} \, dK_{v} \,\Big|\, \F_t \bigg)
\ede
and (\ref{xrt2}) becomes
\bde 
\wh \UV_{u,t,T} = (L_t)^{-1} \, \EP (  L_T \UV_{u,T} \,|\, \F_t ) .
\ede
\el

\proof
The derivation of (\ref{rt1}) is rather standard. Note that the hypothesis~$(H\!P)$ is not needed here and
we may take $t \in [0,T]$. It suffices to take $\UV_{u,T} = g(u) \I_A$ for a Borel measurable map
$g : \rr_+ \to \rr$ and an event $A \in \F_T$ such that
the random variable $\UV_{\tau ,T}= g(\tau ) \I_A$ is $\P$-integrable.
Using part (i) in Lemma \ref{newnn} and the well-known formula for the conditional expectation with respect to $\G_t$,
we obtain
\beq
\lefteqn{ \I_\seq{\tau >t } \EP ( \UV_{\tau ,T} \,|\, \G_t ) =  \I_\seq{\tau >t }
\frac{\EP \big(   \I_\seq{\tau >t } \UV_{\tau ,T} \,|\, \F_t \big)}{\P (  \tau > t \,|\, \F_t)}
= \I_\seq{\tau >t } (G_t)^{-1} \, \EP \Big( \I_A \, \EP \big( \I_\seq{\tau >t }  g(\tau )  \,|\, \F_T \big)\,\big|\, \F_t \Big)} \\
&&=  \I_\seq{\tau >t }(G_t)^{-1} \, \EP \bigg( \I_A \, \int_{(t,\infty ]} g(v) \, dF_{v,T}  \,\Big|\, \F_t \bigg)
=  \I_\seq{\tau >t }(G_t)^{-1} \, \EP \bigg(  \, \int_{(t,\infty ]} g(v) \I_A \, dF_{v,T}  \,\Big|\, \F_t \bigg)\\
&&=  \I_\seq{\tau >t } (G_t)^{-1} \, \EP \bigg(  \, \int_{(t,\infty ]} \UV_{v,T} \, dF_{v,T}  \,\Big|\, \F_t \bigg)
= \I_\seq{\tau >t} \wt \UV_{t,T}
\eeq
where $\wt \UV_{t,T}$ is given by (\ref{rt1}).
We observe that on the event $\{ \tau \leq t\}$ any $\G_t$-measurable
random variable can be represented by a $\sigma (\tau ) \vee \F_t$-measurable random variable $H_{\tau ,t }$.
Let us take any $t \in [0,T)$. To establish (\ref{xrt2}), we need to evaluate $\cE{\P}{U_{\tau,T}}{\G_t}$ on the event $\seq{\tau \leq t}$.
An application of Lemma \ref{gcond} yields
\begin{align*}
\cE{\P}{\I_\seq{\tau\leq t} \UV_{\tau ,T}}{\G_t} &= \lim_{s\downarrow t}\cE{\P}{\I_\seq{\tau\leq t} \UV_{\tau ,T}}{\sigma (\tau ) \vee \F_s}.
\end{align*}
We first compute the conditional expectation $\cE{\P}{\I_\seq{\tau\leq t} \UV_{\tau ,T}}{\sigma (\tau ) \vee \F_s}$  for $0 \leq t< s < T$.
Recall that the hypothesis $(H\!P)$ means that the equality $F_{u,s} F_{s,T} = F_{s,s} F_{u,T}$ holds for all $0 \leq u < s < T$,
which implies that $F_{s,T}\, dF_{u,s} = F_{s}\, dF_{u,T}$ for any fixed $s<T$ and all $u \in [0,t]$.

Hence, for any bounded, $\sigma (\tau ) \vee \F_s$-measurable random variable $H_{\tau ,s}$, we obtain
\beq
\lefteqn{ \EP \big( \I_\seq{\tau\leq t} H_{\tau ,s} \UV_{\tau ,T}  \big)
= \EP \Big( \EP \big(  \I_\seq{\tau\leq t} H_{\tau ,s} \UV_{\tau ,T} \,|\, \F_T \big) \Big)
=\EP \bigg( \int_{[0,t]} H_{u,s} \UV_{u,T}  \, dF_{u,T} \bigg)} \\
&&= \EP \bigg( \int_{[0,t]} H_{u,s} (F_{s})^{-1} F_{s,T} \UV_{u,T} \, dF_{u,s} \bigg)
=  \EP \bigg( \int_{[0,t]} H_{u,s} (F_{s})^{-1} \EP ( F_{s,T} \UV_{u,T}\,|\, \F_s) \, dF_{u,s} \bigg) \\
&&=  \EP \bigg( \int_{[0,t]} H_{u,s} \wh U_{u,s ,T} \, dF_{u,s} \bigg) =  \EP \Big( \I_\seq{\tau\leq t} H_{\tau ,s}  \wh \UV_{\tau ,s ,T} \Big),
\eeq
since $\{\tau \leq t\} \subset \{\tau \leq s \} \subset \{ F_{s} >0\}$, $\P$-a.s.  (see part (ii) in Lemma \ref{newnn}).
This in turn yields
\begin{align*}
\cE{\P}{\I_\seq{\tau\leq t} \UV_{\tau ,T}}{\G_t} &= \lim_{s\downarrow t}\, \cE{\P}{\I_\seq{\tau\leq t} \UV_{\tau ,T}}{\sigma (\tau ) \vee \F_s}
= \lim_{s\downarrow t}\, \I_\seq{\tau\leq t}\wh U_{\tau ,s ,T} \\
						 & = \lim_{s\downarrow t}\,\I_\seq{\tau\leq t}(F_s)^{-1}\, \cE{\P}{F_{s,T}\UV_{u,T}}{\F_s}_{u = \tau }\\
						 & = \I_\seq{\tau\leq t}(F_t)^{-1}\, \cE{\P}{F_{t,T}\UV_{u,T}}{\F_t}_{u = \tau }
= \I_\seq{\tau\leq t}\wh \UV_{\tau ,t, T}
\end{align*}
where the penultimate equality holds by the right-continuity of the filtration $\ff$ and the right-continuity of
processes $F$ and $F_{\cdot ,T}$. This completes the proof of part (i).
For part (ii), we observe that if, in addition, the random field $F_{u,t}$ is completely separable then the asserted formulae
follow from equations (\ref{rt1}) and (\ref{xrt2}).
\endproof

\brem
Let us observe that for $t=0$, we first obtain, on the event $\{ \tau = 0 \}$,
\bde
\EP \big( \I_\seq{\tau = 0 } H_{\tau ,s} \UV_{\tau ,T}  \big)
= \EP \Big( \I_\seq{\tau = 0} H_{\tau ,s}  \wh U_{\tau ,s ,T} \Big)= \EP \Big( \I_\seq{\tau = 0} H_{0,s}  \wh U_{0,s ,T} \Big)
\ede
where, on the event $\{ \tau = 0 \}\subset \{\tau \leq s \} \subset \{ F_s >0\}$,
\bde
\wh \UV_{0,s,T} = (F_s)^{-1} \, \EP (  F_{s,T} \UV_{0,T} \,|\, \F_s ).
\ede
In the second step, we get, on the event $\{ \tau = 0 \} \subset \{ F_0 >0\}$,
\bde
\cE{\P}{\I_\seq{\tau =0} \UV_{\tau ,T}}{\G_0}  =  \lim_{s\downarrow 0}\, \I_\seq{\tau =0}\wh U_{0,s ,T}
=  \I_\seq{\tau = 0 }\wh \UV_{0,0, T}
\ede
where (see (\ref{xrt2}))
\bde
\wh \UV_{0,0,T} = (F_0)^{-1} \, \EP (  F_{0,T} \UV_{0,T} \,|\, \F_0 )
= \frac{\EP ( \P (\tau= 0 \,|\, \F_T) \UV_{0,T} \,|\, \F_0 )}{\P (\tau= 0 \,|\, \F_0)}
\ede
where $F_{0,T} = \P (\tau= 0 \,|\, \F_T)$. We conclude that
\begin{align*}
\EP (  \UV_{\tau ,T} \,|\, \G_0 ) &=  \I_\seq{\tau >0} \frac{1}{\P (\tau >  0 \,|\, \F_0)}  \, \EP \bigg(  \int_{(0,\infty ]} \UV_{u,T} \, d\P (\tau \leq u \,|\, \F_T) \,\Big|\, \F_0 \bigg)
\\ & + \I_\seq{\tau =0} \frac{1}{\P (\tau = 0 \,|\, \F_0)}  \, \EP \bigg(  \int_{[0]} \UV_{u,T} \, d\P (\tau \leq u \,|\, \F_T) \,\Big|\, \F_0 \bigg).
\end{align*}
\erem

\brem For $t=T$, we have that
\be \label{ghgrt}
\EP (  \UV_{\tau ,T} \,|\, \G_T ) =  \I_\seq{\tau > T} \wt \UV_{T,T} + \I_\seq{\tau\leq T}  \UV_{\tau ,T},
\ee
where formula (\ref{rt1}) in Lemma \ref{lmu1} yields
\be \label{crt1}
\wt \UV_{T,T} =  \frac{\EP \big(   \I_\seq{\tau >T } \UV_{\tau ,T} \,\big|\, \F_T \big)}{\P (  \tau > T \,|\, \F_T)}
= (G_T)^{-1} \, \int_{(T,\infty ]} \UV_{u,T} \, dF_{u,T}.
\ee
\erem

\brem
Assume that $F_{u,t}$ satisfies the hypothesis $(H)$. Then formula (\ref{xrt2}) simplifies as follows
\bde
\wh \UV_{u,t,T} = (1-G_t)^{-1} \, \EP (  F_{t,T} \UV_{u,T} \,|\, \F_t )= (F_t)^{-1} \, \EP (  F_{t,t} \UV_{u,T} \,|\, \F_t ) = \EP ( \UV_{u,T} \,|\, \F_t ).
\ede
In particular, if $\wh U_{\tau ,T} = g(\tau)$, we have that $\wh \UV_{u,t,T} = g(u)$ and
\bde
\wt \UV_{t,T} = (G_t)^{-1} \, \EP \big(   \I_\seq{\tau >t} g(\tau ) \,|\, \F_t \big)
= (G_t)^{-1} \, \EP \bigg(  \int_{(t,\infty ]} g(u) \, dF_{u} \,\Big|\, \F_t \bigg).
\ede
\erem

\ssc{Conditional Expectations for Pseudo-Initial Times}

Under the extended density hypothesis, we establish the following counterpart of Lemma \ref{lmu1}.

\bl \label{den1}
Let $\UV_{\cdot ,T} : \rr_+ \times \Omega \to \rr $ be a ${\cal B}(\rr_+) \otimes \F_T$-measurable map.
If $\tau $ is a pseudo-initial time and the random variable $\UV_{\tau ,T}$ is $\P$-integrable
then, for every $t \in [0,T)$,
\bde 
\EP (  \UV_{\tau ,T} \,|\, \G_t ) =  \I_\seq{\tau >t} \wt \UV_{t,T} + \I_\seq{\tau\leq t}\wh \UV_{\tau ,t, T}
\ede
where
\bde 
\wt \UV_{t,T} = (G_t)^{-1} \, \EP \big(   \I_\seq{\tau >t} \UV_{\tau ,T} \,|\, \F_t \big)
= (G_t)^{-1} \, \EP \bigg(  \int_{(t,\infty ]} \UV_{v,T} \, dF_{v,T} \,\Big|\, \F_t \bigg)
\ede
and, for every $0 \leq u \leq t < T$,
\be \label{xrt3}
\wh \UV_{u,t,T} = (\mm_{u,t})^{-1} \, \EP (  \mm_{u,T} \UV_{u,T} \,|\, \F_t ).
\ee
\el

\begin{proof}
It suffices to revise the proof of Lemma \ref{lmu1} on the event $\seq{\tau \leq t}$.
Let us first observe that, by Definition \ref{densityr}, for every $u\in \rr_+$, the process  $(\mm_{u,t})_{t \geq u}$
is a positive $(\P ,\ff)$-martingale and thus  $\{ m_{u,t}=0\} \subset \{ m_{u,s}=0\} \subset \{ m_{u,T} =0 \}$
for all $u \leq t\leq s\leq T$. Recall also that, by convention, we set $0/0 = 0$ and thus $\wh \UV_{u,s,T}$ is well defined.
Therefore, for all $t\leq s\leq T$ and any bounded, $\sigma (\tau ) \vee \F_s$-measurable random variable $H_{\tau ,s}$, we obtain
\beq
\lefteqn{ \EP \big( \I_\seq{\tau\leq t} H_{\tau ,s} \UV_{\tau ,T}  \big)
= \EP \Big( \EP \big(  \I_\seq{\tau\leq t} H_{\tau ,s} \UV_{\tau ,T} \,|\, \F_T \big) \Big)
=\EP \bigg( \int_{[0,\infty ]} \I_\seq{u \leq t} H_{u,s} \UV_{u,T}  \, dF_{u,T} \bigg)} \\
&&= \EP \bigg( \int_{[0,t]} H_{u,s} \mm_{u,T}\UV_{u,T} \, dD_{u} \bigg)
=  \EP \bigg( \int_{[0,t]} H_{u,s} (\mm_{u,s})^{-1}\EP ( \mm_{u,T}\UV_{u,T} \,|\, \F_s) \mm_{u,s}\, dD_{u} \bigg) \\
&&=  \EP \bigg( \int_{[0,t]} H_{u,s} \wh U_{u,s ,T} \, dF_{u,s} \bigg) =  \EP \Big( \I_\seq{\tau\leq t} H_{\tau ,s}  \wh U_{\tau ,s ,T} \Big).
\eeq
By taking limit, and using similar arguments as in the proof of Lemma \ref{lmu1}, we get
\begin{align*}
\cE{\P}{\I_\seq{\tau\leq t} \UV_{\tau ,T}}{\G_t} &= \lim_{s\downarrow t}\, \cE{\P}{\I_\seq{\tau\leq t} \UV_{\tau ,T}}{\sigma (\tau ) \vee \F_s}
= \lim_{s\downarrow t}\, \I_\seq{\tau\leq t}\wh U_{\tau ,s ,T} \\
						 & = \lim_{s\downarrow t}\, \I_\seq{\tau\leq t} (\mm_{u,s})^{-1}\, \EP ( \mm_{u,T}\UV_{u,T} \,|\, \F_s)_{u = \tau}\\
						 & = \I_\seq{\tau\leq t}(\mm_{u,t})^{-1}\, \EP ( \mm_{u,T}\UV_{u,T} \,|\, \F_t)_{u = \tau}
= \I_\seq{\tau\leq t}\wh \UV_{\tau ,t, T},
\end{align*}
which proves \eqref{xrt3}.
\end{proof}

\section{Properties of $\gg$-Local Martingales} \label{sec4}

We consider the map $\wh \UV : \rr^2_+ \times \Omega \to \rr $ and we use the notation
$(u,t,\omega ) \mapsto \wh \UV_{u,t} (\omega )$. We say that $\wh \UV$ is an {\it $\ff$-optional map} when it
is ${\cal B}(\rr_+) \otimes  {\cal O}(\ff)$-measurable, where ${\cal O}(\ff)$ is the $\ff$-optional
$\sigma$-field in $\rr_+ \times \Omega $. In that case, the map $\wh \UV_{\cdot , t}$ is
${\cal B}(\rr_+) \otimes  \F_t $-measurable and the process
$(\wh U_{t,t})_{t \geq 0}$ is $\ff$-optional, in the usual sense.
 We will sometimes need an additional assumption that the process
$(\wh U_{t,t})_{t \geq 0}$ is $\ff$-predictable.

\ssc{$\gg$-Local Martingales for Pseudo-Honest Times}

Let us consider an arbitrary random time $\tau $ such that the process $G=1-F$ is the Az\'ema supermartingale of $\tau $.
We denote by $G = M- \AA $ the Doob-Meyer decomposition of the supermartingale  $G$.
Then the dual $(\P, \ff)$-predictable projection $H^p$ of the indicator process $H_t = \I_{\{ \tau \leq t\}}$ satisfies the equality
 $H^p = \AA$.

The following result, which corresponds to Propositions 5.1 and 5.6 in El Karoui et al. \cite{EJJ}, is an important
step towards establishing a $(\P,\gg)$-semimartingale decomposition of a $(\P,\ff)$-local martingale.

\bt \lab{piu1}
Assume that $\tau $ is a pseudo-honest time and $0 <  F_{u,t} \leq 1$ for every $0< u\leq t$.  Let $\whMM$ be a
$\gg$-adapted and $\P$-integrable process given by the following expression
\be \lab{zz1x}
\whMM_t  =\I_\seq{\tau >t} \wt \UV_t  +\I_\seq{\tau\leq t} \wh \UV_{\tau ,t}
\ee
where  $\wt \UV$ is an $\ff$-adapted, $\P$-integrable process and $\wh \UV$ is an $\ff$-optional map such that for every $t \in \rr_+$
the random variable $\wh \UV_{\tau ,t}$ is $\P$-integrable and the process
$(\wh \UV_t := \wh \UV_{t,t})_{t \geq 0}$ is $\ff$-predictable. Assume, in addition, that the following conditions are satisfied: \hfill \break
(i) the process $(W_{t})_{t \geq 0}$ is a $(\P ,\ff)$-local martingale where
\be \lab{zz2x}
W_t = \wt \UV_t G_t + \int_{(0,t]} \wh \UV_{v} \, dF_v ,
\ee
(ii) for any fixed $u, s \geq 0$, the process $(F_{s,t}\wh \UV^0_{u,t})_{t \geq u\vee s }$
 is a $(\P ,\ff)$-local martingale where we denote $\wh \UV^0_{u,t}= \wh \UV_{u,t} - \wh U_{u,u}$ for every $0 \leq u \leq t$. \hfill \break 
Then the process $(\whMM_t)_{t\geq 0}$ is a $(\P,\gg)$-local martingale.
\et

\begin{proof}
Since the proof proceeds along the similar lines as the proofs of Propositions 5.1 and 5.6 in El Karoui et al. \cite{EJJ},
we will focus on computations and for the details regarding suitable localization and measurability arguments we refer to \cite{EJJ}.
We start by noting that the following decomposition is valid
\begin{align*}
\whMM_t  &= \whMM_t\I_\seq{\tau >t } + \whMM_\tau\I_\seq{\tau \leq t} + (\whMM_t - \whMM_\tau) \I_\seq{\tau \leq t }
 \\ &= \wt \UV_t\I_\seq{\tau >t } + \wh \UV_{\tau,\tau}\I_\seq{\tau \leq t} + (\wh \UV_{\tau,t} - \wh \UV_{\tau,\tau}) \I_\seq{\tau \leq t }.
\end{align*}
It is thus enough to examine the following two subcases, corresponding to conditions (i) and (ii), respectively:
\hfill \break (a) the case of a process $\whMM$ stopped at $\tau $,
\hfill \break (b) the case of a process $\whMM$ such that $\whMM_{\tau \wedge t}=0$ for all $t \geq 0$.

\noindent {\it Case (a).} We first assume that a $\gg$-adapted process $\whMM$ is stopped at $\tau $, specifically,
\be \lab{baru}
 \whMM_t = \I_\seq{\tau >t} \wt \UV_t +\I_\seq{\tau\leq t}  \wh \UV_{\tau }
\ee
where $(\wt \UV_t)_{t \geq 0}$ is an $\ff$-adapted process and $(\wh \UV_t := \wh \UV_{t,t})_{t \geq 0}$ is an $\ff$-predictable process.

\newpage

We start by observing that, for every $0 \leq s <t$,
\begin{align*}
\EP (  \I_\seq{\tau\leq t} \wh \UV_{\tau }  \,|\, \G_s ) &=
\EP (  \I_\seq{s< \tau\leq t} \wh \UV_{\tau }  \,|\, \G_s ) + \EP (  \I_\seq{\tau\leq s} \wh \UV_{\tau }  \,|\, \G_s )\\
&= \I_\seq{\tau > s } (G_s)^{-1} \,\EP \bigg(  \int_{(s,t]} \wh \UV_{v} \, dH_{v} \,\Big|\, \F_s \bigg)+\I_\seq{\tau\leq s}  \wh \UV_{\tau } \\
&= \I_\seq{\tau > s } (G_s)^{-1} \, \EP \bigg(  \int_{(s,t]} \wh \UV_{v} \, d\AA_{v} \,\Big|\, \F_s \bigg)+\I_\seq{\tau\leq s}  \wh \UV_{\tau }\\
&= \I_\seq{\tau > s } (G_s)^{-1} \, \EP \bigg(  \int_{(s,t]} \wh \UV_{v} \, dF_{v} \,\Big|\, \F_s \bigg)+\I_\seq{\tau\leq s}  \wh \UV_{\tau }.
\end{align*}
Therefore, for every $0 \leq s <t$,
\begin{align*}
\EP ( \whMM_t  \,|\, \G_s ) &=  \I_\seq{\tau > s } (G_s)^{-1} \, \EP (  \wt \UV_t G_t  \,|\, \F_s )
+ \I_\seq{\tau > s } (G_s)^{-1} \, \EP \bigg(  \int_{(s,t]} \wh \UV_{v} \, dF_{v} \,\Big|\, \F_s \bigg) + \I_\seq{\tau\leq s}  \wh \UV_{\tau } \\
&=  \I_\seq{\tau > s } (G_s)^{-1} \, \EP (  W_t - W_s  \,|\, \F_s ) + \I_\seq{\tau > s }
(G_s)^{-1} \, \EP (  \wt \UV_s G_s  \,|\, \F_s ) + \I_\seq{\tau\leq s}  \wh \UV_{\tau }
\\ &= \I_\seq{\tau >s} \wt \UV_s +\I_\seq{\tau\leq s}  \wh \UV_{\tau }  =  \whMM_s
\end{align*}
where we used the assumption that the process $W$ given by (\ref{zz2x}) is a $(\P ,\ff)$-martingale.
We conclude that the process  $\whMM$  given by (\ref{baru}) is a $(\P ,\gg)$-martingale.

\noindent {\it Case (b).} Let us denote $\wh \UV^0_{v,t}= \wh \UV_{v,t} - \wh \UV_{v,v}$ for $0 \leq v \leq t$.
Consider a $\gg$-adapted process $ \whMM$ given by $\whMM_t = \I_\seq{\tau\leq t}\wh \UV^0_{\tau ,t} $ where $\wh \UV^0_{t,t}=0$.
We need to show that the equality $\EP (  \I_\seq{\tau\leq t} \wh \UV^0_{\tau ,t} \,|\, \G_s ) = \I_\seq{\tau\leq s} \wh \UV^0_{\tau ,s} $
holds for every $0 \leq s < t$. From part (i) in Lemma \ref{lmu1}, we obtain, for every $0\leq  s < t$,
\begin{align*}
\EP (  \I_\seq{\tau\leq t} \wh \UV^0_{\tau ,t}  \,|\, \G_s ) & = \I_\seq{\tau > s } (G_s)^{-1} \,
\EP \bigg(  \int_{(s,t]} \wh \UV^0_{v,t} \, dF_{v,t} \,\Big|\, \F_s \bigg)
\\ &+ \I_\seq{\tau\leq s} (F_s)^{-1} \, \EP (  F_{s,t} \wh \UV^0_{u,t} \,|\, \F_s )_{| u = \tau } = I_1 + I_2 .
\end{align*}
Let us examine $I_1$. We first assume that $s>0$. Recall that we assume that the hypothesis $(H\!P)$ holds and
$0 <  F_{v,t} \leq 1$ for every $0< v \leq t$.
Hence, for $0<s \leq v \leq t$, we can write $dF_{v,t} =  F_{s,t}\, d(F_{v,v}F_{s,v}^{-1}) = F_{s,t}\, dD^s_v$ where the process
$(D^s_v = F_{v,v}F_{s,v}^{-1})_{v \geq  s}$ is increasing and $\ff$-adapted. Consequently,
\begin{align*}
I_1 &= \I_\seq{\tau > s } (G_s)^{-1} \, \EP \bigg( \int_{(s,t]} \EP ( F_{s,t}\wh \UV^0_{v,t}\, |\, \F_v ) \, dD^s_v \,\Big|\, \F_s \bigg)\\
&= \I_\seq{\tau > s } (G_s)^{-1} \, \EP \bigg(\int_{(s,t]} F_{s,v}\wh \UV^0_{v,v} \, dD^s_{v} \,\Big|\, \F_s \bigg) = 0
\end{align*}
where we first used condition (ii) and subsequently the equality $\wh \UV^0_{v,v}=0$.
It remains to examine the case $s = 0$. We denote $\wh \UV^{0+}_{\tau,t} = \max \, (\wh \UV^0_{\tau,t},0)$
and $\wh \UV^{0-}_{\tau,t} = \max \, ( - \wh \UV^0_{\tau,t},0)$. Then,
for all $t >0$,
\begin{align*}
I_1 & = \I_\seq{\tau > 0 } (G_0)^{-1} \,\EP \bigg(  \int_{(0,t]} \wh \UV^0_{v,t} \, dF_{v,t} \,\Big|\, \F_0 \bigg)
= \I_\seq{\tau > 0 } (G_0)^{-1} \,\EP \bigg(  \lim_{s \downarrow 0} \cE{\P}{\I_\seq{s \leq \tau \leq t}\wh \UV^0_{\tau,t}}{F_t}\,\Big|\, \F_0 \bigg) \\
    & = \I_\seq{\tau > 0 } (G_0)^{-1} \,\EP \bigg(  \lim_{s\downarrow 0} \cE{\P}{\I_\seq{s\leq \tau \leq t}\wh \UV^{0+}_{\tau,t}}{F_t}\,\Big|\, \F_0 \bigg)
     - \EP \bigg(  \lim_{\epsilon \downarrow 0} \cE{\P}{\I_\seq{s \leq \tau \leq t}\wh \UV^{0-}_{\tau,t}}{F_t}\,\Big|\, \F_0 \bigg).
\end{align*}
By the monotone convergence theorem for conditional expectations, we obtain
\begin{align*}
I_1 & = \I_\seq{\tau > 0 } (G_0)^{-1} \,\lim_{s \downarrow 0}\EP \bigg(  \int_{(s,t]} \wh \UV^0_{v,t} \, dF_{v,t} \,\Big|\, \F_0 \bigg) \\
    & = \I_\seq{\tau > 0} (G_0)^{-1} \, \lim_{s\downarrow 0}\EP \bigg( \int_{(s,t]} \EP ( F_{s,t}\wh \UV^0_{v,t}\, |\, \F_u ) \, dD^s_v \,\Big|\, \F_0 \bigg) = 0
\end{align*}
where we used condition (ii) and the equality $\wh \UV^0_{v,v}=0$.
For $I_2$, using again condition (ii), we obtain, for $0 \leq  u \leq s < t$,
\bde
I_2 = \I_\seq{\tau\leq s} (F_s)^{-1} \, \EP (  F_{s,t} \wh \UV^0_{u,t} \,|\, \F_s )_{| u = \tau }  =
\I_\seq{\tau\leq s} (F_s)^{-1} \, F_{s,s} (\wh \UV^0_{u,s})_{| u = \tau } = \I_\seq{\tau\leq s} \wh \UV^0_{\tau ,s}.
\ede
We conclude that the process $\big(\I_\seq{\tau\leq t} \wh \UV^0_{\tau ,t}\big)_{t \geq 0}$ is a $(\P,\gg)$-martingale
and thus the proof of the proposition is completed.
\end{proof}

The following corollary to Theorem \ref{piu1} deals with the special case when the process $\UV$ given by (\ref{zz1x})
is continuous at $\tau $. It is easy to check that under the assumptions of Corollary \ref{remi} the process
$(\wh U_t := \wh U_{t,t})_{t \geq 0}$ is $\ff$-predictable.

\bcor \label{remi}
Under the assumptions of Theorem \ref{piu1} we postulate, in addition,
that the equality $\wt U_{t-} = \wh U_{t,t}$ holds for every $t \in \rr_+$. Then the process $\whMM$ is continuous at $\tau$ and
condition $(i)$ in Theorem \ref{piu1} can be replaced by the following condition: \hfill\break
$(i')$ the process $(W_{t})_{t \geq 0}$ is a $(\P,\ff)$-local martingale where
\be
W_t = \wt \UV_t G_t + \int_{(0,t]} \wt \UV_{u-} \, dF_u .
\ee
\ecor

To establish another corollary to Theorem \ref{piu1}, we assume that $F_{u,t}$ is completely separable.

\bcor \lab{xeed}
Under the assumptions of Theorem \ref{piu1} we postulate, in addition, that the $(\P,\ff)$-conditional distribution of
$\tau$ satisfies $F_{u,t} = K_uL_t$ for all $0\leq u \leq t$, where $K$ is a positive, $\ff$-adapted,
increasing process and $L$ is a positive $(\P,\ff)$-martingale. Then
condition (ii) in Theorem \ref{piu1} can be replaced by the following condition: \hfill\break
(ii') For every $u\geq 0$, the process $(W_{u,t}= L_t \wh \UV^0_{u,t})_{t\geq u}$ is a $(\P,\ff)$-martingale.
\ecor

\begin{proof}  We observe that the computations for case (b) in the proof of Proposition  \ref{piu1}
can be simplified. Using part (ii) in Lemma \ref{lmu1}, we obtain, for every $0 \leq s < t$ (note that $L_t=0$ on the event $\{L_s=0\}$)
\beq \lefteqn{
\EP (  \I_\seq{\tau\leq t} \wh \UV^0_{\tau ,t}  \,|\, \G_s ) = \I_\seq{\tau > s } (G_s)^{-1} \,
\EP \bigg(  \int_{(s,t]} L_t \wh \UV^0_{u,t} \, dK_{u} \,\Big|\, \F_s \bigg)
+ \I_\seq{\tau\leq s} (L_s)^{-1} \, \EP (  L_t \wh \UV^0_{u,t} \,|\, \F_s )_{| u = \tau } } \\
&&= \I_\seq{\tau > s } (G_s)^{-1} \, \EP \bigg(  \int_{(s,t]} \EP ( W_{u,t}\, |\, \F_u ) \, dK_{u} \,\Big|\, \F_s \bigg)
+\I_\seq{\tau\leq s} (L_s)^{-1} \, \EP (  W_{u,t} \,|\, \F_s )_{| u = \tau } \\
&&= \I_\seq{\tau > s } (G_s)^{-1} \, \EP \bigg(  \int_{(s,t]} W_{u,u} \, dK_{u} \,\Big|\, \F_s \bigg)
+  \I_\seq{\tau\leq s} (L_s)^{-1} (W_{u,s})_{| u = \tau } = \I_\seq{\tau\leq s} \wh \UV^0_{\tau ,s}
\quad \quad \quad \quad \quad \quad \quad \quad \quad \quad \quad \quad \quad \quad \quad \quad
\eeq
where we used condition (ii') in the penultimate equality
and the equality $W_{u,u}= L_u \wh \UV^0_{u,u}=0$ in the last one.
\end{proof}

\ssc{$\gg$-Local Martingales for Pseudo-Initial Times}

It was necessary to assume in Theorem \ref{piu1} that the $(\P,\ff)$-conditional distribution $F_{u,t}$
is non-degenerate, since the random measure $D^s_u:=F_{u,u}(F_{s,u})^{-1}$ is not always
well defined when the $(\P,\ff)$-conditional distribution $F_{u,t}$ is degenerate. In order to circumvent
this technical assumption, one can postulate instead that $F_{u,t}$ satisfies the hypothesis $(E\!D)$.
In the next result, we work under the setup of Theorem \ref{piu1}, but we no longer assume that
$0 < F_{u,t} \leq 1$ for every $0< u\leq t$.

\bp \lab{piu5}
Suppose that $\tau $ is a pseudo-initial time and \eqref{qaz} holds with a positive random field  $(\mm_{s,t})_{t \geq s}$ and an $\ff$-adapted
increasing process $D$. Then condition (ii) in Theorem \ref{piu1} can be replaced by the following condition: \hfill \break
(ii*) for every $u\geq 0$, the process $(\mm_{u,t} \wh \UV^0_{u,t} )_{t \geq u}$ is a $(\P ,\ff)$-local martingale.
\ep

\begin{proof}
We only need to adjust the proof of Theorem \ref{piu1} in case (b). Let $\whMM_t = \I_\seq{\tau\leq t}\wh \UV^0_{\tau ,t}$ where $\wh U^0_{t,t} = 0$.
Using Lemma \ref{den1}, we obtain, for every $0 \leq s < t$ (recall that $\{ m_{u,s}=0\} \subset \{ m_{u,t}=0\}$)
\begin{align*}
\EP (  \I_\seq{\tau\leq t} \wh \UV^0_{\tau ,t}  \,|\, \G_s ) & = \I_\seq{\tau > s } (G_s)^{-1} \,
\EP \bigg(  \int_{(s,t]} \wh \UV^0_{u,t} \, dF_{u,t} \,\Big|\, \F_s \bigg)
\\ &+ \I_\seq{\tau\leq s} \, (\mm_{u,s})^{-1}\EP (  \mm_{u,t} \wh \UV^0_{u,t} \,|\, \F_s )_{| u = \tau } = I_1 + I_2.
\end{align*}
 The integral $I_1$ satisfies
\begin{align*}
I_1 &= \I_\seq{\tau > s } (G_s)^{-1} \, \EP \bigg( \int_{(s,t]} \wh \UV^0_{u,t} \, dF_{u,t} \,\Big|\, \F_s \bigg)
= \I_\seq{\tau > s } (G_s)^{-1} \, \EP \bigg(\int_{(s,t]} \mm_{u,t}\wh \UV^0_{u,t} \, dD_{u} \,\Big|\, \F_s \bigg)\\
&= \I_\seq{\tau > s } (G_s)^{-1} \, \EP \bigg(\int_{(s,t]} \EP \big( \mm_{u,t}\wh \UV^0_{u,t} \,|\, \F_u \big) \, dD_{u} \,\Big|\, \F_s \bigg) = 0
\end{align*}
where to obtain the last equality we first used assumption (ii*) and next the equality $\wh \UV^0_{u,u}=0$. The integral $I_2$ simplifies to
\bde
I_2 = \I_\seq{\tau\leq s} (\mm_{u,s})^{-1} \, \EP (  \mm_{u,t} \wh \UV^0_{u,t} \,|\, \F_s )\big|_{u = \tau }  =
\I_\seq{\tau\leq s} (\mm_{u,s})^{-1}  \, \mm_{u,s} (\wh \UV^0_{u,s})\big|_{u = \tau } = \I_\seq{\tau\leq s} \wh \UV_{\tau ,s}.
\ede
We conclude that the process $\whMM $ is a $(\P,\gg)$-martingale, as was required to show.
\end{proof}

\section{Compensators of the Indicator Process}  \label{sec5}

Our next goal is to compute the $(\P,\gg)$-dual predictable projection of the indicator process $H_t = \I_\seq{\tau \leq t}$
where, as usual, we denote by $\gg$ the progressive enlargement of $\ff$ with a random time $\tau$.
Recall that the Doob-Meyer decomposition of $G$ is denoted by
$G=  M - \AA $. Then the $(\P,\ff)$-dual predictable projection  (i.e., the $(\P,\ff)$-{\it compensator}) of $H$,
denoted as $H^p$, coincides with the $\ff$-predictable,
increasing process $\AA$. To find the $(\P,\gg)$-dual predictable projection (i.e., the $(\P,\gg )$-{\it compensator}) of $H$, it is enough
to apply the following classic result, due to Jeulin and Yor \cite{JY1} (see also Guo and Zeng \cite{GZ}), and to compute explicitly
the $(\P, \ff )$-compensator of~$H$.

\bt
Let $\tau$ be a random time with the Az\'ema supermartingale $G$. Then the $(\P,\gg)$-{\it compensator} of $H$ equals
\be \lab{JYfor}
H^{p,\gg}_t  =  \int_{(0,t\wedge \tau]} \frac{1}{G_{u-}} \, dH^p_u ,
\ee
meaning that the process $H - H^{p,\gg}$ is a $(\P,\gg)$-martingale.
\et

\ssc{Compensator of $H$ under Complete Separability}

Let us first examine the case where the $(\P,\ff)$-conditional distribution of $\tau$ under $\P$
is completely separable, that is, $F_{u,t} = K_uL_t$.
We assume, in addition, that the increasing process $K$ is $\ff$-predictable. It is worth noting
that both assumptions are satisfied in the construction of $\tau $ based on a predictable multiplicative system,
provided that $G_t<1$ for all $t> 0$ (see \cite{LR1}).
By applying the integration by parts formula to $F$ and using the assumption that $K$ is an $\ff$-predictable process,
we obtain
\bde
F_t = K_tL_t = K_0L_0 + \int_{(0,t]} K_u \, dL_u + \int_{(0,t]} L_{u-}\, dK_u .
\ede
Hence, by the uniqueness of the Doob-Meyer decomposition, we conclude that $d\AA_u = L_{u-}dK_u$.
Consequently, using the Jeulin-Yor formula (\ref{JYfor}), we obtain
\bde
H^{p,\gg}_t  = \int_{(0,t\wedge \tau]} \frac{L_{u-}}{1-L_{u-}K_{u-}} \, dK_u .
\ede

\ssc{Compensator of $H$ for a Pseudo-Initial Time}

Assume now that $\tau $ is a pseudo-initial time, that is, the hypothesis $(E\!D)$ holds (see Definition \ref{densityr}).
Let the process $m$ be given as $( m_t := m_{t,t})_{t \geq 0}$ and let $\mp$ denote the $(\P,\ff)$-predictable
projection of $m$. For brevity, we will use the notation $X \stackrel{mart} = Y$ whenever $X-Y$ is a $(\P,\ff)$-local martingale.
In the following, it is assumed that the process $\int m_u \, dD_u$ is of integrable variation.

\bl \label{m}
The process $m^D_t := \int_{(0,t]} (m_{u,t} - m_{u}) \, dD_u$ is a $(\P,\ff)$-martingale.
\el

\begin{proof}
 By splitting the integral and taking the conditional expectation, we obtain, for any $s\leq t$,
\begin{align*}
\cE{\P}{m^D_t}{\F_s} &= \EP \bigg( \int_{(s,t]} (m_{u,t} - m_{u}) \, dD_u \,\Big|\, \F_s \bigg)
 +  \EP \bigg( \int_{(0,s]} (m_{u,t} - m_{u})  \, dD_u \,\Big|\, \F_s \bigg).
 \\ & = \EP \bigg( \int_{(s,t]} \EP (m_{u,t} - m_{u} \,|\, \F_u ) \, dD_u \,\Big|\, \F_s \bigg)
 +  \EP \bigg( F_{s,t} - F_{0,t} -\int_{(0,s]}m_{u} \, dD_u \,\Big|\, \F_s \bigg)
 \\ & = F_{s,s} - F_{0,s} - \int_{(0,s]}m_{u} \, dD_u = \int_{(0,s]}(m_{u,s} - m_{u}) \, dD_u = m^D_s ,
\end{align*}
since, for any fixed $s$, the process $F_{s,t} = \int_{[0,s]} m_{u,t}\, dD_u $ is assumed to be
a $(\P ,\ff)$-martingale (see Definitions \ref{hyp2} and \ref{densityr}). \end{proof}

\bp
 Assume that $\tau $ is a pseudo-initial time and the process $\mm$ is a special semimartingale with the canonical decomposition
$\mm = N + P$ where $N$ is the local martingale part. If the predictable covariation of $N$ and $D-D^p$ exists then
the $(\P,\ff)$-compensator of $H$  is given by the formula
\be \label{dres}
H^p_t = \pbracket{N}{D-D^p}_t+ \int_{(0,t]}  {}^p \mm_{u} \, dD^p_u
\ee
where  ${}^p \mm$ is the $(\P,\ff)$-predictable projection of $\mm$.
\ep

\begin{proof}
It suffices to compute the Doob-Meyer decomposition of $F$. Using (\ref{qaz}) and the canonical decomposition of $\mm$, we obtain
\begin{align*}
F_t & = \int_{[0,t]} m_{u,t}\, dD_u =  \mm_{0,t}\Delta D_0 + \int_{(0,t]} (\mm_{u,t} - \mm_{u,u})\, dD_u + \int_{(0,t]} N_u \,dD_u + \int_{(0,t]} P_u \,dD_u  \\
    & = \mm_{0,t}\Delta D_0 + m^D_t + \sbracket{N}{D}_t + \int_{(0,t]} (N_{u-} + P_u) \,dD_u
\end{align*}
where in the second equality we used the fact that $\sbracket{N}{D}_t = \int_{(0,t]}\Delta N_{u} \, dD_u$
(see Proposition 9.3.7.1 in \cite{JYC}). It is clear that the process $\mm_{0,t}\Delta D_0$ is a $(\P,\ff)$-martingale, whereas $m^D_t$ is a $(\P,\ff)$-martingale in view of Lemma \ref{m}. Using the $(\P,\ff)$-dual predictable projection of $D$, we can write
\begin{align*}
\int_{(0,t]} (N_{u-} + P_u) \,dD_u  & = \int_{(0,t]} (N_{u-} + P_u) \,d(D_u - D^p_u) + \int_{(0,t]} (N_{u-} + P_u) \,dD^p_u \\
& \stackrel{\rm mart} = \int_{(0,t]} (N_{u-} + P_u) \,dD^p_u .
\end{align*}
Finally, the $\ff$-predictable covariation of $N$ and $D-D^p$ is assumed to exist and
\bde
\sbracket{N}{D} = \sbracket{N}{D- D^p} + \sbracket{N}{D^p} \stackrel{\rm mart} = \pbracket{N}{D-D^p}
\ede
where the second equality holds since: \hfill \break
(a) $N$ and $D-D^p$ are $(\P,\ff)$-local martingales so that $\sbracket{N}{D- D^p}-\pbracket{N}{D-D^p}$ is $(\P,\ff)$-local martingale, \hfill \break (b) by Y{\oe}urp's lemma (see the proof of Proposition 9.3.7.1 in \cite{JYC}), the process  $\sbracket{N}{D^p}$ is a $(\P,\ff)$-local martingale since
 $N$ is $(\P,\ff)$-local martingale and $D^p$ is a predictable process of finite variation.

 We have thus shown that
\begin{align*}
F_t  \stackrel{\rm mart} = \pbracket{N}{D-D^p}_t + \int_{(0,t]} (N_{u-} + P_u) \,dD^p_u
\end{align*}
where the right-hand side is the $\ff$-predictable process of finite variation.
Now, let $(T_n)_{n\in \nn}$ be a localizing sequence such that $M$ and all the $(\P,\ff)$-local martingales defined above are
 uniformly integrable $(\P ,\ff)$-martingales once stopped at $T_n$. Then we can conclude that, for every $n \in \nn$,
\begin{align*}
\AA^{T_n}_t &= \pbracket{N}{D-D^p}^{T_n}_t + \int_{(0,t]} (N_{u-} + P_u) \I_{\llb 0, T_n \rrb}(u)\,dD^p_u\\
&= \pbracket{N}{D-D^p}^{T_n}_t + \int_{(0,t]} \,^p m_u\I_{\llb 0, T_n \rrb}(u) \,dD^p_u,
\end{align*}
where the last equality follows from the equalities: ${}^p (\mm^{T_n})  ={}^p \big(N^{T_n}\big)+ {}^p \big( P^{T_n} \big)
= N^{T_n}_{-} + P^{T_n}$ (see, in particular, Theorem 4.5 in  \cite{N}) and (note that the process $\I_{\llb 0,T_n\rrb}$ is
$\ff$-predictable)
\bde
\,^p \big(m^{T_n}\big)\I_{\llb 0,T_n\rrb}= \,^p \big(m^{T_n}\I_{\llb 0,T_n\rrb} \big)= \,^p \big( m\I_{\llb 0,T_n\rrb}\big)
= \,^p m\I_{\llb 0,T_n\rrb}.
\ede
To complete the proof of equality (\ref{dres}), it suffices to let $n \rightarrow \infty$.
\end{proof}

\brem
As an example, consider the case where the process $D$ is predictable. Then we obtain from \eqref{dres}
\bde
A_t = H^p_t =\int_{(0,t]} {}^p \mm_{u}\,dD_u
\ede
and it is enough to require that ${}^p \mm$ exists.
\erem

\section{Hypothesis $(H')$ and Semimartingale Decompositions} \label{sec6}

The aim of this section is to analyze the validity of the classic hypothesis $(H')$ for progressive enlargements
associated with pseudo-honest and pseudo-initial times. We establish here the main results of this work,
Theorems \ref{piu} and \ref{piu7}, and we study the case of the multiplicative construction of a random
time associated with a predetermined Az\'ema submartingale.

Let us first recall the general definition,
in which $\kk$ stands for any enlargement of $\ff$, that is, any filtration such that $\ff \subset \kk$.

\bd
The {\it hypothesis $(H')$} is said to hold for $\ff $ and its enlargement $\kk$ whenever
any $(\P,\ff)$-semimartingale is also a $(\P,\kk )$-semimartingale.
\ed

For exhaustive studies of the hypothesis $(H')$ the interested reader is referred to Jeulin \cite{J2}, who examined a general case as well
as honest times, and Jacod  \cite{JJ}, who worked under the density hypothesis and covered the initial times.
The latter study was recently extended by Kchia and Protter \cite{KP}, who dealt with the progressive
enlargement with a general stochastic process, and not only the indicator process of a random time.

As is well known, to establish the hypothesis $(H')$ between $\ff$ and any enlargement $\kk$,
it suffices to show that any bounded $(\P,\ff)$-martingale is a $(\P,\kk )$-semimartingale (see Yor \cite{Y1}).
This crucial  observation follows, for instance, from the Jacod-M\'emin decomposition of a $(\P,\ff)$-semimartingale; $X=X_0+K+B+N$
where $K$ represents large jumps, $B$ is predictable of finite variation and $N$ is a local
martingale with jumps bounded by 1 (see, e.g., page 3 in Jeulin \cite{J2}).
One can then show that, under the hypothesis $(H')$, any bounded $(\P,\ff)$-martingale
is in fact a special $(\P,\kk)$-semimartingale and this in turn implies that, more generally,
any special $(\P,\ff)$-semimartingale remains a special $(\P,\kk )$-semimartingale.
Therefore, assuming that the hypothesis $(H')$ holds for $\ff$ and $\kk$, the natural goal is thus to find the canonical semimartingale decomposition
with respect to the enlarged filtration $\kk$ of a bounded $(\P,\ff)$-martingale. In addition,
if the hypothesis $(H')$ for $\ff$ and $\kk$ fails to hold then the goal is to describe the class
of $(\P,\ff)$-semimartingales that remain also $(\P,\kk)$-semimartingales.
For an exhaustive study of these problems, the reader may consult Jeulin~\cite{J2}.

Let us now focus on the case of the progressive enlargement of a filtration $\ff$ with a random time.
It is well known, in particular, that for any random time $\tau $ with values in $\err_+$
and any $(\P,\ff)$-local martingale $U$, the stopped process $U^{\tau }$ is a $(\P,\gg)$-semimartingale
(see Yor \cite{Y1} and Jeulin and Yor \cite{JY1}). For the future reference, we first state a theorem, which recalls and combines
results from papers by Jeulin and Yor \cite{JY1} (for part (i), see Theorem 1 in \cite{JY1}) and
Jeulin \cite{J2} (for part (ii), see Proposition 4.16 in \cite{J2}).

\bt \lab{projy}
Let $\gg$ be the progressive enlargement of $\ff$ with
an arbitrary random time $\tau $. If $U$ is a $(\P,\ff )$-local martingale
then the stopped process $U^{\tau }$ is a $(\P,\gg)$-special semimartingale. \hfill \break
\noindent (i)  The process
\be \lab{dsa1}
U_{t\wedge \tau} - \int_{(0,t \wedge \tau ]}\frac{1}{G_{u-}}\, d( \langle U, M \rangle_u+ \breve U^p_u )
\ee
is a $(\P, \gg )$-local martingale, where $\breve  U^p$ stands for the dual $\ff$-predictable projection of the process
$\breve  U_t = \Delta U_\tau\I_\seq{\tau \leq t}$.  \hfill \break
(ii)  The process
\be \lab{dsa2a}
U_{t\wedge \tau} - \int_{(0,t \wedge \tau ]}\frac{1}{G_{u-}}\, d\langle U, \bar M \rangle_u
\ee
is a $(\P, \gg )$-local martingale, where $\bar M$ is the unique BMO martingale such that
$\EP (N_{\tau }) = \EP ( N_{\infty } \bar M_{\infty })$ for every bounded $(\P,\ff)$-martingale $N$.
\et

\begin{proof}
We shall only demonstrate part (ii), since the proof of part (i) is readily available in \cite{JY1}.
Let us suppose, without loss of generality, that $U$ is a $(\P,\ff)$-martingale such that $U_0 = 0$.
It is well known that for every $\gg$-stopping time $\wh T$ there exists an $\ff$-stopping time $T$ such
that $\tau\wedge \wh T = \tau\wedge T$. Therefore,
\bde
 \expect{\P}{U_{\wh T\wedge \tau}}  = \expect{\P}{U_{T\wedge \tau}}
 = \expect{\P}{\int_{[0,\infty]} U^T_{u} d\bar A_u} = \expect{\P}{U^T_{\infty}\bar M_\infty}
 = \expect{\P}{\int_{(0,\infty]} \frac{\I_\seq{\tau \geq u}}{G_{u-}} \, d\pbracket{U^T}{\bar M}_u}
\ede
where the last equality follows from the integration by parts formula and the property of the dual $\ff$-predictable projection
(in addition, it is known that ${}^p(\I_\seq{\tau \geq u})= G_{u-}$; see, e.g., page 576 in Jeulin \cite{J1}).  This in turn implies
\bde
 \expect{\P}{U_{\wh T\wedge \tau} - \int_{(0,\tau\wedge \wh T]} \frac{1}{G_{u-}} \, d\pbracket{U}{\bar M}_u} = 0 = \expect{\P}{U_{0\wedge \tau}}
\ede
and thus the process defined in \eqref{dsa2a} is a $(\P,\gg)$-local martingale.
\end{proof}

Observe that \eqref{dsa1} and \eqref{dsa2a} yield alternative representations for the canonical decomposition
of the special $(\P,\gg)$-semimartingale $U^{\tau}$. Hence, from the uniqueness of the canonical decomposition of
a special semimartingale, we deduce the equality
\be \label{dsa4}
\int_{(0,t \wedge \tau ]}\frac{1}{G_{u-}}\, d( \pbracket{U}{M}_u + \breve U^p_u )
 = \int_{(0,t \wedge \tau ]}\frac{1}{G_{u-}}\, d\pbracket{U}{\bar M}_u
\ee
which necessarily holds for an arbitrary $(\P,\ff )$-local martingale $U$.

\brem \lab{yy6}
It is known that $G = \bar M - \bar \AA$ where $\bar \AA = H^o$ is the dual $\ff$-optional projection of $H$.
Under the assumption $({\bf C})$ (that is, when all $\ff$-martingales are continuous) and/or
the assumption $({\bf A})$ (that is, when the random time $\tau $ avoids all $\ff$-stopping times),
we have that $\bar M = M$ and thus also $H^p = H^o$ or, equivalently, $\AA = \bar \AA$.
\erem

We will also need the following auxiliary result (see Lemma 4 in Jeulin and Yor \cite{JY1}).

\bl  \lab{corxc}
Let $U$ be a $(\P,\ff)$-local martingale. Denote by $\breve U^p$ the dual $(\P, \ff)$-predictable projection of the process
\bde
\breve U_t = \int_{(0,t]} \Delta U_s \, dH_s = \Delta U_\tau\I_\seq{\tau \leq t}.
\ede
Then the process
\be \lab{pcoi1}
\int_{(0,t]} \Delta U_s \, dH_s- \int_{(0,t \wedge \tau ]} \frac{1}{G_{u-}} \, d\breve U^p_u=
\Delta U_\tau\I_\seq{\tau \leq t}  - \int_{(0,t \wedge \tau ]} \frac{1}{G_{u-}} \, d\breve U^p_u
\ee
is a $(\P, \gg) $-local martingale.
\el

\ssc{Hypothesis $(H')$ for the Progressive Enlargement}

Our goal is to show that if a random time $\tau $ given on some filtered probability space $(\Omega, \F, \ff,\P)$
satisfies the hypothesis $(H\!P)$ and its $(\P,\ff)$-conditional distribution is positive
 then the $\gg$-semimartingale decomposition of a $(\P,\ff)$-semimartingale can be computed explicitly.
From Theorem  \ref{projy}, we know that for any $(\P,\ff )$-local martingale $U$, the stopped process $U^{\tau }$ is a $(\P,\gg)$-special
semimartingale with explicitly known canonical decomposition. Therefore, it will be enough to focus on the behavior of a process $U$ after $\tau $.
It should be stressed that we do not claim that the hypothesis $(H\!P)$ implies that the hypothesis $(H')$
between $\ff$ and the progressive enlargement
$\gg$ holds, in general, since certain additional assumptions will be imposed when deriving alternative
versions of $\gg$-semimartingale decompositions of a  $(\P,\ff )$-local martingale.
 Theorem \ref{piu} furnishes an explicit $(\P,\gg)$-semimartingale decomposition of a $(\P ,\ff)$-local martingale $U$
when $\tau $ is a pseudo-honest time. This result can be seen as a counterpart
of Proposition 5.9 in El Karoui et al. \cite{EJJ} who dealt with the case of an initial time (see also Jeanblanc and Le Cam
\cite{JL} who examined both initial and honest times). In the special case of the multiplicative
construction (see Corollary \ref{munn1} and Remark \ref{jean}), it is also related to Theorem 7.1 in Jeanblanc and Song \cite{JS1}.

\bt \lab{piu}
Assume that $\tau $ is a pseudo-honest time such that, for every $s \geq 0$, the bounded $(\P,\ff)$-martingale $(F_{s,u})_{u \geq s}$
is strictly positive. If $U$ is a $(\P ,\ff)$-local martingale then the process $\whMM$ is a $(\P,\gg )$-local martingale where
\be \label{gdecomp}
\whMM_t = U_t - \int_{(0,t \wedge \tau ]} (G_{u-})^{-1}\, d( \langle U, M \rangle_u + \breve U^p_u )
 - \int_{(t \wedge s ,t]}  (F_{s,u-})^{-1} \, d\langle U,F_{s,\cdot}\rangle_u  \Big|_{s = \tau} .
\ee
Hence $U$ is a $(\P,\gg)$-special semimartingale and equality \eqref{gdecomp} yields its canonical decomposition.
\et

Note that $(t \wedge \tau ,t] = \emptyset $  on the event $\{ \tau \geq t \}$, since manifestly
$(t \wedge s ,t] = \emptyset $ for all $s \geq t$.
Before proceeding to the proof of Theorem \ref{piu}, we make some pertinent remarks and prove a preliminary lemma.

\brem Recall that any local martingale is locally in the space ${\cal H}^1$. If $N$ is a BMO martingale then,
by Fefferman's inequality (see Revuz and Yor \cite{RY}), there exists a constant $c$ such that for any local martingale $U$
\bde
 \EP \left( \int_0^{\infty } | d [U,N]_t | \right) \leq c \, \| U \|_{{\cal H}^1} \| N \|_{BMO}.
\ede
Consequently, the process $[U,N]$ is locally of integrable variation and its compensator $\langle U,N \rangle $ is well defined.
\erem

\brem The Az\'ema's supermartingale $G$ is generated by the $\ff$-predictable, increasing process $\AA $,
in the sense that, for every $t\geq 0$,
\bde
G_t = \cE{\P}{\AA_{\infty}}{\F_t} - \AA_t = M_t - \AA_t.
\ede
This implies that the process $M_t = \cE{\P}{\AA_\infty}{\F_t}$ is a BMO martingale since $G \leq 1$ (see Proposition 10.13 in \cite{HWY}).
It is also known that any bounded martingale is a BMO martingale (see Proposition 10.11 in \cite{HWY}).
\erem

As a first step towards the proof of Theorem \ref{piu}, we  establish
the existence of the integrals appearing in right-hand side of \eqref{gdecomp}.

\bl \lab{vvff}
Under the assumptions of Theorem \ref{piu}, the integrals in the right-hand side of equality
\eqref{cgdecomp} are well defined.
\el

\begin{proof} The process $(G_t)^{-1} \I_\seq{\tau > t}$ is known to have, with probability 1, finite
left-hand limits for all $t \in \rr_+$ (see Yor \cite{Y1}) and thus it possess a finite left-hand limit at $\tau $. Hence the first integral
in \eqref{gdecomp} is a well defined $\gg$-adapted process of finite variation.
Next, let us check that for all $u\leq t$ the integral $Z_{u,t} = \int_{(u,t]} (F_{u,v-})^{-1}d\pbracket{U}{F_{u,\cdot}}_v$ is well defined.
We proceed as follows. Under the standing assumption that $(F_{u,t})_{t \geq u}$ is a strictly positive process,
the stochastic logarithm
\bde
\mathcal{L}(F_{u,\cdot}) = \int_{(u,t]} (F_{u,v-})^{-1} dF_{u,v}
\ede
is well defined. For the existence of the predictable bracket $\pbracket{U}{\mathcal{L}(F_{u,\cdot})}$, it is sufficient to check that the c\`adl\`ag process
$\mathcal{L}(F_{u,\cdot})$ is a locally bounded martingale and for this purpose it is enough to show that
the jump process $\Delta \mathcal{L}(F_{u,\cdot})_t = (F_{u,t-})^{-1}\Delta F_{u,t}$
is locally bounded for $t \geq u$. The latter property is clear since the left-continuous process
$(F_{u,t-})^{-1}$ is locally bounded (we use here the property that $(F_{u,t})_{t \geq u}$
is a strictly positive $(\P,\ff)$-martingale) and the jumps of $F_{u,\cdot }$ are obviously bounded by~1.
We conclude that the integral $Z_{u,t}$ is well defined since
\bde
Z_{u,t}= \int_{(u,t]} (F_{u,v-})^{-1}d\langle U , F_{u,\cdot}\rangle_v = \int_{(u,t]}d\langle U , \mathcal{L}(F_{u,\cdot})\rangle_v .
\ede
Moreover, the process $(Z_{u,t})_{t\geq u}$ is of locally integrable variation, since
\bde
\expect{\P}{\int_u^{\infty } \, |d\langle U , \mathcal{L}(F_{u,\cdot})\rangle_v |}
\leq c \, \| U \|_{{\cal H}^1} \| \mathcal{L}(F_{u,\cdot}) \|_{BMO},
\ede
where the local martingale $U$ is locally in ${\cal H}^1$ and the locally bounded martingale $\mathcal{L}(F_{u,\cdot})$ is locally in the space BMO.  \end{proof}

We are now in a position to prove Theorem \ref{piu}.  Note that Theorem \ref{projy} is not employed in the proof
of Theorem \ref{piu}, although we use Lemma \ref{vvff} to compensate the jump of $U$ at $\tau $.

\noindent {\it Proof of Theorem \ref{piu}.} Although the present set-up is more general than the one
studied in El Karoui et al. \cite{EJJ}, some steps in our proof
are analogous to those employed in the proof of Proposition 5.9 in  \cite{EJJ}. We start by noting that
we may and do assume, without loss of generality, that $\MM$ is a uniformly integrable $(\P ,\ff)$-martingale.
In view of Lemma \ref{corxc}, it suffices to show that if $U$ is a $(\P ,\ff)$-local martingale then the process $\whMM$
is a $(\P,\gg )$-local martingale where
\be \label{cgdecomp}
\whMM_t = U_t - \int_{(0,t \wedge \tau ]} (G_{u-})^{-1} \, d\langle U, M \rangle_u - \Delta U_\tau\I_\seq{\tau \leq t}
 - \int_{(t \wedge s ,t]}  (F_{s,u-})^{-1} \, d\langle U,F_{s,\cdot}\rangle_u  \Big|_{s = \tau} .
\ee
We note that $\whMM$ satisfies the equality $\whMM = \MM - \BB^*$ where the process $\BB^*$ is defined by
\bde
\BB^*_t  = \wt \BB_t \I_\seq{\tau >t} +  \wh \BB_{\tau ,t} \I_\seq{\tau\leq t}
\ede
where in turn the $\ff$-predictable process $\wt \BB$ equals
\be \label{piu3}
\wt \BB_t = \int_{(0,t]}(G_{u-})^{-1}\, d\langle U, M \rangle_u
\ee
and the map $\wh \BB$ is given by, for all $0 \leq u \leq t$,
\bde
\wh \BB_{u,t} = \wt \BB_{u}+\Delta U_u + \int_{(u,t]} (F_{u,v-})^{-1} \, d\langle \MM , F_{u,\cdot }\rangle_v
\ede
where $\Delta U_u = U_u - U_{u-}$.

\newpage
We need to show that the process
\bde
\whMM_t =  \I_\seq{\tau >t} \wt \UV_t +\I_\seq{\tau\leq t}  \wh \UV_{\tau ,t }
= (\MM_t - \wt \BB_t ) \I_\seq{\tau >t} + (\MM_t - \wh \BB_{\tau ,t}) \I_\seq{\tau\leq t}
\ede
is a $(\P, \gg)$-local martingale. We observe that $\wh \UV_{u,t}= \MM_t - \wh \BB_{u,t}$ is an $\ff$-optional map
and the process $\wh U_{t,t} = U_{t-} - \wt \BB_t$ is $\ff$-predictable. Hence the assumptions
of Theorem \ref{piu1} are satisfied.

Therefore,  in order to show that $\whMM$ is a $(\P, \gg)$-local martingale,
it suffices to show that: (i) the process
\be \lab{waq}
(U_t-\wt \BB_t)G_t - \int_{(0,t]}(U_{u-}-\wt \BB_{u}) \, dG_u
\ee
is a $(\P,\ff)$-local martingale and (ii) for any fixed $u,s \geq 0$, the process
$(F_{s,t}\wh U^0_{u,t})_{t \geq s\vee u}$ is a $(\P,\ff)$-local martingale.
For brevity, we will use the notation $X \stackrel{mart} = Y$ whenever the process $X-Y$ is a $(\P,\ff)$-local martingale.

To establish (i), we observe that
\beq
d\big((U_t-\wt \BB_t)G_t\big) &=& U_{t-}\, dG_t + G_{t-}\, dU_t + d[U,G]_t
 - \wt \BB_{t}\, dG_t - G_{t-} \, d \wt \BB_t  \\
 &\stackrel{mart} =& U_{t-}\, dG_t + d[U,G]_t- \wt \BB_{t}\, dG_t - G_{t-} \, d \wt \BB_t
\eeq
Consequently,
\bde
d\big((U_t-\wt \BB_t)G_t\big) - (U_{t-}-\wt \BB_{t-})\, dG_t  \stackrel{mart} =
 d[U,G]_t  - d\langle U , M \rangle_u =  d[U,B]_t + d[U,M]_t  - d\langle U , M \rangle_u .
\ede
Hence the process given by (\ref{waq}) is a $(\P,\ff)$-local martingale
since, by Y{\oe}urp's lemma, the process $[U,B]_t$ is a $(\P,\ff)$-local martingale
and the process $[U,M]  - \pbracket{U}{M}$ is  a $(\P,\ff)$-local martingale as well. We conclude that
the property (i) is valid in the present setup.

For (ii), we fix $u \geq 0$ and we observe that, for every $t \geq u$,
\be \lab{nexw}
\wh U^0_{u,t} = U_t - \wh \BB_{u,t} - \wh U_{u,u} =  U_{t} - Z_{u,t} - \VV_u
\ee
where we denote $\VV_u =   \wt \BB_{u} + \Delta U_u + \wh U_{u,u}$
and  $(Z_{u,t})_{t\geq u}$ is the $\ff$-predictable process given by
\be  \lab{nexw1}
Z_{u,t}  = \int_{(u,t]} (F_{u,v-})^{-1} \, d\langle \MM , F_{u,\cdot }\rangle_v .
\ee
We fix $u,s \geq 0$. By applying the integration by parts formula, we obtain, for $t \geq s$,
\bde
d( F_{s,t} \wh U^0_{u,t} ) =  U_{t-} \, dF_{s,t} + F_{s,t-}\, dU_{t} +
d[\MM,F_{s,\cdot}]_t - Z_{u,t} \, dF_{s,t} - F_{s,t-}\, dZ_{u,t} - \VV_u \, dF_{s,t}
\ede
where we used the fact that the process $(Z_{u,t})_{t\geq u}$ is $\ff$-predictable.
Recall that $U$ and $(F_{s,t})_{t\geq s}$ are $(\P,\ff)$-martingales.
Therefore, to show that the process $( F_{s,t} \wh U^0_{u,t} )_{t \geq u \vee s}$ is a $(\P,\ff)$-local martingale, it is
enough to check that, for all $u \leq s \leq t$ and $s \leq u \leq t$,
\be
d[\MM,F_{s,\cdot}]_t - F_{s,t-}\, dZ_{u,t} \stackrel{\rm mart} =
d \pbracket{\MM }{F_{s,\cdot}}_t - F_{s,t-}\, dZ_{u,t}
\stackrel{\rm mart} =  0 . \label{piu2}
\ee
Since we assumed that $\tau $ is a pseudo-honest time, using Remark \ref{drivingM}, we obtain, for all $u\leq s \leq t$,
\bde
F_{s,t-}\, dZ_{u,t} = \frac{F_{s,t-}}{F_{u,t-}} \, d \pbracket{\MM}{F_{u,\cdot }}_t = \frac{F_{s,t-}}{F_{u,t-}} \,
d\pbracket{ \MM}{\frac{F_{u,s}}{F_{s,s}}\, F_{s,\cdot}}_t =
\frac{F_{s,t-}F_{u,s}}{F_{u,t-}F_{s,s}} \, \pbracket{\MM}{F_{s,\cdot }}_t = d\pbracket{\MM}{F_{s,\cdot }}_t .
\ede
Similarly, for all $s\leq u \leq t$, we get 
\bde
F_{s,t-}\, dZ_{u,t} = \frac{F_{s,t-}}{F_{u,t-}} \, d \pbracket{\MM}{F_{u,\cdot }}_t =
\frac{F_{s,u}F_{u,t-}}{F_{u,u}F_{u,t-}} \,  d \pbracket{\MM}{F_{u,\cdot }}_t 
= d\pbracket{\MM}{\frac{F_{s,u}}{F_{u,u}} \, F_{u,\cdot}}_t   = d\pbracket{\MM}{F_{s,\cdot }}_t .
\ede
This shows that $d\pbracket{\MM}{F_{s,\cdot }}_t - F_{s,t-}\, dZ_{u,t} =0$ and thus the second equality \eqref{piu2} is trivially satisfied.
This means, of course, that the property (ii) is satisfied.
Using Theorem \ref{piu1}, we thus conclude that the process $\whMM$ given by
\eqref{cgdecomp} is a $(\P,\gg )$-local martingale. This in turn implies that
the process $\whMM$ given by \eqref{gdecomp} is a $(\P,\gg )$-local martingale, as was required to show.
\hfill $\Box$

The next result is borrowed from Kchia et al. \cite{KLP},  who examined the case of any two enlargements that coincide
after a random time $\tau $. However, for simplicity of presentation, their result (see Theorem 3 in \cite{KLP}) is stated here for the
special case of the progressive enlargement $\gg$ and the initial enlargement $\gg^*$, which are known to
coincide after $\tau $. For brevity, we write hereafter $\wt C = \langle U,M \rangle + \breve U^p$ where $\breve U^p$
is the dual $\ff$-predictable projection of the process $\breve U_t = \Delta U_\tau\I_\seq{\tau \leq t}$.

\bt \lab{kchia}
Let $U$ be a $(\P, \ff )$-local martingale. Suppose that $B$ is a $\gg^*$-predictable process of finite variation such that $U-B$ is a
$(\P,\gg^*)$-local martingale. Then the process
\be
U_t - \int_{(0,t \wedge \tau ]}\frac{1}{G_{u-}}\, d\wt C_t
 - \int_{(t \wedge \tau ,t]} dB_u
\ee
is a $(\P, \gg )$-local martingale.  Hence $U$ is a $(\P,\gg)$-special semimartingale.
\et

\ssc{Special Cases of Pseudo-Honest Times} \lab{sec62}

From Theorem \ref{piu}, we know that the hypothesis $(H')$ holds for a pseudo-honest time with a strictly positive
$(\P,\ff)$-conditional distribution. Moreover, this result furnishes also a general expression for the canonical decomposition
with respect to $\gg$ of an arbitrary $(\P,\ff)$-local martingale. Of course, any $\ff$-predictable process
of finite variation is also a $\gg$-predictable process of finite variation and thus it suffices to focus
on the canonical decomposition with respect to $\gg$ of a $(\P,\ff)$-local martingale, rather than a
$(\P,\ff)$-special semimartingale.
Our next goal is to examine some useful consequences of Theorem \ref{piu} under alternative additional assumptions imposed
on a pseudo-honest time under consideration. In particular, we will compare various semimartingale decompositions for pseudo-honest times
with their classic counterparts established for honest times by Barlow \cite{B2}, Jeulin and Yor \cite{JY1}, and Jeulin and Yor \cite{JY2}.

For the readers's convenience, we first recall the most pertinent results regarding the case of an honest time.
It was shown by Barlow \cite{B2} (Theorem 3.10) and Yor \cite{Y1} (Theorem 4) that the hypothesis $(H')$ holds for $\ff$ and
its progressive enlargement with an honest time.
The following result summarizes the well known properties of the progressive enlargement for an honest
time (see Theorem A in Barlow \cite{B2}, Theorem 2 in Jeulin and Yor \cite{JY1}, and
Theorem 15 in Jeulin and Yor \cite{JY2}).

\bt \lab{thhon}
Let $\gg$ be the progressive enlargement of $\ff$ with an honest time $\tau $. If $U$ is a $(\P,\ff )$-local martingale
then $U$ is a $(\P,\gg)$-special semimartingale. \hfill \break
\noindent (i) The process
\be \lab{dsa3}
U_t - \int_{(0,t \wedge \tau ]} (G_{u-})^{-1} 
\, d\wt C_u + \int_{(t \wedge \tau ,t]} (F_{u-})^{-1}\, d\wt C_u ,
\ee
is a $(\P, \gg)$-local martingale. \hfill \break
(ii) The process
\be \lab{dsa3cc}
 U_t - \int_{(0,t \wedge \tau ]} (G_{u-})^{-1} \, d\langle U,\bar M \rangle_u + \int_{(t \wedge \tau ,t]}
 (F_{u-})^{-1}\, d\langle U,\bar M \rangle_u
\ee
is a $(\P, \gg)$-local martingale, where $\bar M$ is the $(\P,\ff)$-martingale of class BMO introduced in part (ii) of Theorem \ref{projy}.
\et

Equalities \eqref{dsa3} and \eqref{dsa3cc} yield alternative representations for the canonical decomposition of
$U$ with respect to the progressive enlargement $\gg$.
 Therefore, for an honest time $\tau$, we obtain the following equality complementing \eqref{dsa4}
\be\label{aftertau}
\int_{(t \wedge \tau ,t]} (F_{u-})^{-1}\, d\wt C_u = \int_{(t \wedge \tau ,t]} (F_{u-})^{-1}\, d\langle U,\bar M \rangle_u
\ee
which holds for any $(\P,\ff )$-local martingale $U$.

\sssc{Completely Separable Case}

As a first special case of a pseudo-honest time, we consider the situation where the $(\P,\ff)$-conditional
distribution of a random time $\tau$ is completely separable (see Definition \ref{defsep}).
Then we obtain the following immediate corollary to Theorem \ref{piu}. Corollary \ref{piu4} will later
be exemplified through the predictable multiplicative construction of a random time (see Corollary \ref{munn1}).

\bcor \label{piu4}
Under the assumptions of Theorem \ref{piu}, we postulate, in addition, that
the $(\P,\ff)$-conditional distribution of $\tau$ is completely separable,
that is, $F_{u,t} = K_uL_t$ for every $0\leq u\leq t$ where $L$ is a strictly positive $(\P,\ff)$-martingale.
If $\MM$ is a $(\P,\ff)$-local martingale then the process $\whMM$ is a $(\P,\gg )$-local martingale where
\be \lab{iiuu}
\whMM_t =  U_t - \int_{(0,t \wedge \tau ]} (G_{u-})^{-1} \, d\wt C_u 
- \int_{(t \wedge \tau ,t]} (L_{u-})^{-1} \, d\langle U,L\rangle_u .
\ee
\ecor

\begin{proof}
We note that under the assumptions of Corollary \ref{piu4}, formula \eqref{cgdecomp}
reduces to \eqref{iiuu}.
\end{proof}

\sssc{Predictable Multiplicative Construction}

The next corollary shows that if $\tau $ is a pseudo-honest time with a non-degenerate $(\P,\ff)$-conditional
distribution then, under certain technical assumptions, the $(\P,\gg)$-semimartingale decomposition of a $(\P,\ff)$-local martingale is
analogous to the one derived by other authors for an honest time and reported in Theorem \ref{thhon}.
It is worth stressing that the present setup is manifestly different from the one covered by
Theorem \ref{thhon} and, in fact, our results do not cover the case of an honest time.

Let us observe that, under the assumptions of Theorem \ref{piu},
for every $s \geq 0$, the process $(C_{s,u} = (F_u)^{-1}F_{s,u})_{u \geq s}$ is positive, decreasing, and
$\ff$-adapted. Indeed, from \eqref{hpx} we obtain the equality
$(F_u)^{-1}F_{s,u} = (F_{u,t})^{-1}F_{s,t}$, which holds for all $0 \leq s < u <t$, where $(F_{u,t})_{u \geq 0}$ is an increasing process.

\bcor \label{hphonest}
Let the assumptions of Theorem \ref{piu} be satisfied. Assume, in addition, that, for every $s \geq 0$, the decreasing process
$(C_{s,u} = (F_u)^{-1}F_{s,u})_{u \geq s}$ is $\ff$-predictable. If $U$ is a $(\P,\ff)$-local martingale then the process $\whMM$ given by
\be \label{vgdecomp}
\whMM_t = U_t - \int_{(0,t \wedge \tau ]}  (G_{u-})^{-1}\, d\wt C_u
 +\int_{(t \wedge \tau ,t]} ({}^pF_u)^{-1}\, d\langle U,M \rangle_u
\ee
is a $(\P,\gg)$-local martingale.
\ecor

\begin{proof} To show that the second integral in the right-hand side of \eqref{vgdecomp}
is a well-defined process of locally integrable variation, we observe that $0< F_{u-} \leq {}^pF_{u}$
(since $F$ is a submartingale) and thus
\begin{align*}
&\expect{\P}{\int_0^{\infty } \I_\seq{\tau < u } ({}^pF_{u})^{-1} \,d|\langle \MM,M \rangle|_t } =
\expect{\P}{\int_0^{\infty } F_{u-}({}^pF_{u})^{-1} \,| d\langle \MM,M \rangle_t|} \\
& \leq \expect{\P}{\int_0^{\infty } \,| d\langle \MM,M \rangle_t |}
 \leq c \, \| U \|_{{\cal H}^1} \| M \|_{BMO}
\end{align*}
where we used Fefferman's inequality in the last inequality.
To obtain \eqref{vgdecomp} from  \eqref{gdecomp}, we start by noting that $ C_{s,u}$ is an $\ff$-predictable multiplicative system associated
with a positive submartingale $F$ (see Meyer \cite{M}) or \cite{LR1}).
By assumption, the process $F_t$ is strictly positive for $t>0$ and thus, by Theorem 4.1 and Corollary 4.1 in \cite{LR1},
the unique $\ff$-predictable multiplicative system $C_{s,u}$ associated with $F$ satisfies
the following stochastic differential equations
\be \label{xoo2}
dC_{s,u} =  - C_{s,u-}(\Fp_u)^{-1}d\AA_u = -C_{s,u}(F_{u-})^{-1}d\AA_u
\ee
where the second equality follows from the equality $C_{s,u}\Fp_u =  C_{s,u-}F_{u-}$ (see formula (17) in \cite{LR1}).
Since the decreasing process $(C_{s,u})_{u \geq s}$ is assumed to be $\ff$-predictable, the integration by parts formula yields,
for any fixed $s \geq 0$,
\bde
dF_{s,u} = C_{s,u}\, dF_u + F_{u-}\, dC_{s,u} = - C_{s,u} \, dM_u
\ede
where the second equality follows from \eqref{xoo2} and the Doob-Meyer decomposition $dF_t = dA_t-dM_t$.
We only need to focus on the last term in formula \eqref{gdecomp}. We have, for all $s \leq t$
\bde
\int_{(s,t]} (F_{s,{u-}})^{-1} \, d\langle U,F_{s,\cdot}\rangle_u
 = -  \int_{(s,t]}  (F_{u-}C_{s,u-})^{-1}  C_{s,u}\, d\langle U,M \rangle_u
= - \int_{(s,t]} ({}^pF_u)^{-1}\, d\langle U,M \rangle_u ,
\ede
as was required to show. \end{proof}

\brem
Let us consider the situation of Corollary \ref{munn1} and let us assume, in addition, that the avoidance property
({\bf A}) holds, that is, $\tau $ avoids all $\ff$-stopping times. Then ${}^p F = F_{-}$ and the
$(\P,\ff)$-local martingale $U$ is continuous at $\tau $.
Hence  \eqref{vgdecomp} becomes
\be \lab{nvfr4}
 \whMM_t = \MM_t - \int_{(0,t \wedge \tau ]} (G_{u-})^{-1} \, d\langle \MM,M \rangle_u
 +  \int_{(t \wedge \tau,t]} (F_{u-})^{-1} \,d\langle \MM,M \rangle_u
\ee
since, under assumption ({\bf A}), the martingales $M$ and $\bar M$ are known to coincide as well.
Note that under assumption ({\bf A}) alternative semimartingale decompositions \eqref{dsa3} and \eqref{dsa3cc} obtained for an honest time
reduce to \eqref{vgdecomp} as well.
\erem

\brem
For every finite $\gg$-stopping time $T$, we obtain
\begin{align*}
\int_{(0,T]}\I_\seq{\tau < u}({}^pF_{u})^{-1} \,d\langle \MM,M \rangle_t & = \int_{(0,T]}\I_\seq{\tau < u}(F_{u-})^{-1} \,d\langle \MM,M \rangle_t  + \int_{(0,T]} \frac{\I_\seq{\tau < u}\Delta A_u}{{}^pF_{u}F_{u-}} \,d\langle \MM,M \rangle_t \\
& = \int_{(0,T]}\I_\seq{\tau < u}(F_{u-})^{-1} \,d\langle \MM,M \rangle_t  + \int_{(0,T]} \frac{{}^p(\I_{\llb \tau \rrb})_u\I_\seq{\tau < u}}{{}^p(\I_{\llb \tau ,\infty \rrb})_uF_{u-}} \,d\langle \MM,M \rangle_t .
\end{align*}
We focus now on the second integral
\begin{align*}
\int_{(0,T]} \frac{{}^p(\I_{\llb \tau \rrb}\I_{\llb \tau,\infty \rrb})_u\I_\seq{\tau < u}}{{}^p(\I_{\llb \tau ,\infty \rrb})_uF_{u-}} \,d\langle \MM,M \rangle_u & = \int_{(0,T]} {}^p\left(\frac{\I_{\llb \tau \rrb}}{\I_{\llb \tau ,\infty \rrb}}\right)_u\frac{\I_\seq{\tau < u}}{F_{u-}} \,d\langle \MM,M \rangle_u
\end{align*}
where the first equality holds since for any bounded measurable process $X$ and $Y$, we have ${}^p(XY^{-1}) = {}^p(XY)({}^pY)^{-1}$,
with the convention that $0/0 = 0$ (see \cite{J1}). Therefore,  another representation of the integral after $\tau$ is given by
\begin{align*}
\int_{(0,t]}\I_\seq{\tau < u}({}^pF_{u})^{-1} \,d\langle \MM,M \rangle_t
& = \int_{(0,t]}\I_\seq{\tau < u}(F_{u-})^{-1} \,d\langle \MM,M \rangle_t
+ \int_{(0,t]} {}^p\left(\frac{\I_{\llb \tau \rrb}}{\I_{\llb \tau ,\infty \rrb}}\right)_u
\frac{\I_\seq{\tau < u}}{F_{u-}} \,d\langle \MM,M \rangle_u .
\end{align*}
\erem

We will now describe a particular instance where the assumptions of Corollary \ref{munn1} are satisfied.
For the reader's convenience, we will first summarize the main steps in an explicit construction of
a random time associated with an arbitrary Az\'ema submartingale $F$,
as developed in \cite{LR1}. We now assume that we are given a predetermined Az\'ema submartingale $F$, that is,
an arbitrary submartingale $F= (F_t)_{t \in \err_+}$, defined on a filtered probability space $(\Omega, \F , \ff,\P)$,
satisfying the inequalities $0 \leq F_t \leq 1$ for every $t \in \rr_+$ and with $F_\infty =1$.
The predictable multiplicative construction of a random time $\tau $ associated with $F$ runs as follows:
\begin{itemize}
\item We start by establishing the existence of an $\ff$-predictable multiplicative
system $\wh C_{u,t}$ associated with a positive submartingale $F$ (see Meyer \cite{M} and Theorem 4.1 in \cite{LR1}).
\item Subsequently, using a (possibly non-unique) $\ff$-predictable multiplicative
system $\wh C_{u,t}$, we define the unique $(\P,\ff)$-conditional distribution $\wh F_{u,t}$ by setting
(see Theorem 4.2 and Lemma 5.1 in \cite{LR1})
\bde 
\wh F_{u,t} =
\begin{cases}
\cE{\P}{F_u}{\F_t}, & t \in  [0,u),\\
\wh C_{u,t} F_t, & t \in [u,\infty].
\end{cases}
\ede
\item Finally, we construct a random time $\tau $ on the extended probability space
$(\wh{\Omega}, \wh \F , \wh{\ff}, \wh\P)$
such that $\wh \P (\tau \leq t\,|\, \F_t )=F_t$ for all $t \in \rr_+$ (see Theorem 5.1 in \cite{LR1}).
\end{itemize}

Recall that $\wh \P$ is chosen in such a way that the probability measures $\wh \P$ and $\P$ coincide on the filtration $\ff$.
 It is also worth noting that if $G=M-A$ is the Doob-Meyer decomposition of $G$ then, using
the uniqueness of the Doob-Meyer decomposition, we deduce that $M_t = \EPH ( H^p_{\infty } \,|\, \F_t )$ and $A =H^p$.
It is clear that Corollary \ref{hphonest} can now be applied to the random time $\tau $ constructed as above.
Specifically, we are in a position to establish the following result, in which we use the
fact that, under stronger assumptions on $F$, the unique $\ff$-predictable multiplicative
system $\wh C_{u,t}$ associated with $F$ is known explicitly (see \cite{LR1}).

\bcor \label{munn1}
Assume that $F_t >0 $ and $F_{t-} >0$ for every $t >0$ and a pseudo-honest time $\tau $
is constructed using the unique $\ff$-predictable multiplicative system associated with $F$.
If $\MM$ is a $(\wh \P ,\ff)$-local martingale then the process $\whMM$ given by
\eqref{vgdecomp} is a $(\wh \P,\gg)$-local martingale.
\ecor

\begin{proof}
The statement follows immediately from Corollary \ref{hphonest}. Alternatively, it can also be deduced
from Corollary \ref{piu4}. To see this, we start by noting that Proposition 5.1 in \cite{LR1} implies
that $F_{u,t}$ is completely separable with $L = F \mathcal{E}$ where $\mathcal{E}$ is the Dol\'eans exponential (see formula (30) in \cite{LR1})
\bde
\mathcal{E}_t =  \mathcal{E}_t \bigg(-\int_{(0,\, \cdot \, ]} (\Fp_s)^{-1} \, d\AA_s \bigg),
\ede
so that $d\mathcal{E}_{t} =- \mathcal{E}_{t-} (\Fp_t)^{-1} \, d\AA_t$.
It is also known that  $ \mathcal{E}_{t}\Fp_t = \mathcal{E}_{t-}F_{t-} = L_{t-}$
(this is a consequence of formula (17) in \cite{LR1}). Since  $\mathcal{E}$ is an $\ff$-predictable process of finite
variation, by applying the integration by parts formula, we obtain
\bde
dL_t=  {\mathcal{E}}_{t} \, dF_t + F_{t-} \, d\mathcal{E}_t
 =  (\Fp_t)^{-1} L_{t-} \, dF_t - F_{t-} \mathcal{E}_{t-} (\Fp_t)^{-1} \, d\AA_t .
\ede
Consequently, we also have that
\begin{align*}
(L_{t-})^{-1}\, dL_t  = -(\Fp_t)^{-1}\,dM_t + (\Fp_t)^{-1}\,d\AA_t
+ (L_{t-})^{-1} F_{t-} \mathcal{E}_{t-} (\Fp_t)^{-1} \, d\AA_t
				 = -(\Fp_t)^{-1}\,dM_t
\end{align*}
and this in turn implies
\bde
(L_{t-})^{-1}\, d\langle U,L \rangle_t  = -({}^pF_t)^{-1}\, d\langle U,M \rangle_t.
\ede
To conclude the proof, it suffices to apply Corollary \ref{piu4}.
\end{proof}

\brem \lab{jean}
Corollary \ref{munn1} corresponds to Theorem 7.1 in Jeanblanc and Song \cite{JS1} who work under the
assumption that $G_t = N_t e^{-\Lambda_t}$ where $N$ is a positive local martingale and $\Lambda $
 is a continuous increasing process. Consequently, the martingale
part $M$ in the Doob-Meyer decomposition of $G$ satisfies $dM_t = e^{-\Lambda_t } \, dN_t$.
Moreover, they postulate that that $G_t <1 $ and $G_{t-} < 1 $ for every $t >0$.
It is shown in \cite{JS1} that one may construct a random time $\tau$ on the product space $\Omega \times \err_+$
 and with respect to a suitably defined probability measure $\Q$
such that the equality $\Q=\P$ is satisfied on $\ff$ and $\Q (\tau >t\,|\, \F_t)=G_t$
for all $t \in \err_+$. Jeanblanc and Song also show  (see Theorem 7.1 in \cite{JS1}) that if a process $U$
is a $(\P,\ff)$-local martingale then the process $\whMM$ given by \eqref{nvfr4}
is a $(\P,\gg)$-local martingale. More precisely, under their assumptions, formula \eqref{nvfr4} becomes
\bde 
 \whMM_t = \MM_t - \int_{(0,t \wedge \tau ]} \frac{e^{-\Lambda_u}}{G_{u-}} \, d\langle \MM, N \rangle_u
 +  \int_{(t \wedge \tau,t]} \frac{e^{-\Lambda_u}}{F_{u-}} \,d\langle \MM, N \rangle_u .
\ede
For a more general result in this vein, see also Theorems 4.2 and 4.4 in  Jeanblanc and Song \cite{JS2}.
It is worthwhile to observe that, in the setup considered in \cite{JS1}, the canonical solution $\tau$ satisfies
the `local density hypothesis in canonical form' or in another words,  the  hypothesis (ED) is satisfied
(see Theorem 5.1 in \cite{JS1}). This implies
that for any $\ff$-stopping time $T$ (for the definition of ${\cal E}_t(u)$, see Corollary 3.1 in \cite{JS1})
\bde
\Q (\tau = T\,|\, \F_\infty ) = \int_{\llb T\rrb} N_u \mathcal{E}_\infty(u) e^{-\Lambda_u}\, d\Lambda_u
 =N_T\mathcal{E}_\infty(T) e^{-\Lambda_T}\Delta\Lambda_T = 0,
\ede
where we also used the assumption made in \cite{JS1} that the process $\Lambda$ is continuous. We conclude that
the avoidance property ({\bf A}) holds. As a consequence, the dual $\ff$-predictable and $\ff$-optional projections of
$H$ coincide and thus also $M = \bar M$.
In our general setting, the assumption that $\Lambda$ is continuous is equivalent
to the assumption that the $\ff$-predictable process $\AA $ which generates $G$ is continuous.
Under this assumption, Theorem 7.1 in \cite{JS1} can be recovered from Corollary \ref{munn1}.
\erem

\sssc{Optional Multiplicative Construction}

In the next special case, we assume that a pseudo-honest time $\tau $ is constructed using an $\ff$-optional multiplicative system associated with
a predetermined Az\'ema submartingale $F$.
We suppose that we are given an Az\'ema submartingale $F$. In order to construct
an $\ff$-optional multiplicative system associated with $F$, we proceed as in
 Section 4.3 in \cite{LR1}. Specifically, we start by assuming that $F$
is given as $F_t = \P ( \widehat \tau \leq t \,|\, \F_t)$ for
some random time $\widehat \tau $. Next, an $\ff$-optional multiplicative system associated with $F$ is defined.
We may assume, without loss of generality, that an auxiliary random time $\widehat \tau $ was constructed using an
$\ff$-predictable multiplicative system $(\wh C_{u,t})_{u,t\geq 0}$ associated with $F$. Hence this additional requirement
is not restrictive.

More formally, to construct an $\ff$-optional multiplicative system associated with $F$,
we set $\wh H = \I_{\llb \wh \tau , \infty \llb}$ and we define $\wh A = \wh H^o$
and $\wh M_t = \EP ( \wh H^o_{\infty } \,|\, \F_t )$. As in Remark \ref{yy6},
we note that the equality $G = \wh M - \wh A$ yields an $\ff$-optional decomposition of $G$.
We define the random field $(C_{u,t})_{u,t \in \err_+}$ by setting $C_{u,t}=1$ for all $u \geq t$
and, for all $t \geq u$,
\be \label{oo2}
dC_{u,t} =  - C_{u,t-}(F_t)^{-1}\, d\wh \AA_t .
\ee
Then, from Corollary 4.2 in \cite{LR1}, the random field $(C_{s,u})_{s,u \in \err_+}$ is an $\ff$-optional multiplicative system
associated with $F$. The $(\P,\ff)$-conditional distribution of a random time is now defined by
\bde 
F_{u,t} =
\begin{cases}
\cE{\P}{F_u}{\F_t}, & t \in  [0,u),\\
 C_{u,t} F_t, & t \in [u,\infty].
\end{cases}
\ede
Finally, the random time $\tau $ can be constructed using once again Theorem 5.1 in \cite{LR1}. It is then not difficult
to check that $\tau $ is a pseudo-honest time. It is important to emphasize that
we do not claim that the equality  $\wh \AA  = H^o$ holds where, as usual, we write $H = \I_{\llb \tau , \infty \llb}$.
Therefore, the $\ff$-optional decomposition of $G= \bar M - \bar A$, which is obtained for a random time $\tau $ as outlined in Remark \ref{yy6},
does not coincide with  the $\ff$-optional decomposition $G = \wh M - \wh A$, which is associated with an auxiliary random time $\wh \tau $.

\bcor \label{nhphonest}
Let the assumptions of Theorem \ref{piu} be satisfied by a pseudo-honest time $\tau $ constructed using
the $\ff$-optional multiplicative system given by \eqref{oo2}.
If $U$ is a $(\P,\ff)$-local martingale then the process
\be \label{nvgdecomp}
U_t - \int_{(0,t \wedge \tau ]} (G_{u-})^{-1}\, d\langle U, \bar M \rangle_u
 + \int_{(t \wedge \tau ,t]} (F_{u-})^{-1}\, d\langle U, \wh M \rangle_u
\ee
is a $(\P,\gg)$-local martingale.
\ecor

\begin{proof} The first integral in \eqref{nvgdecomp} is dealt with as in part (ii) of Theorem  \ref{projy}.
For the second integral in \eqref{nvgdecomp}, we start by noting that it
is a well-defined process of locally integrable variation, since
\begin{align*}
& \expect{\P}{\int_0^{\infty } \I_\seq{\tau < t } (F_{t-})^{-1} \,d|\langle \MM,\wh M \rangle|_t } =
\expect{\P}{\int_0^{\infty } F_{t-}(F_{t-})^{-1} \,d|\langle \MM,\wh M \rangle|_t } \\
& \leq \expect{\P}{\int_0^{\infty } \,d|\langle \MM, \wh M \rangle|_t } \leq c \, \| U \|_{{\cal H}^1} \| \wh M \|_{BMO}
\end{align*}
where the local martingale $U$ is locally in ${\cal H}^1$ and $\wh M$ is the BMO martingale.
Since the process $(C_{s,u})_{u \geq s}$ is decreasing, the integration by parts formula yields,
for any fixed $s \geq 0$,
\bde
dF_{s,u} = C_{s,u-}\, dF_u + F_{u}\, dC_{s,u} = -  C_{s,u-} \, d\wh M_u
\ede
where the second equality follows from \eqref{oo2} and the decomposition  $G =\wh M - \wh \AA$.
Hence, for all $s \leq t$,
\bde
\int_{(s,t]} (F_{s,{u-}})^{-1} \, d\langle U,F_{s,\cdot}\rangle_u
 = -  \int_{(s,t]}  (F_{u-}C_{s,u-})^{-1}  C_{s,u-}\, d\langle U,\wh M \rangle_u
= - \int_{(s,t]} (F_{u-})^{-1}\, d\langle U,\wh M \rangle_u .
\ede
To conclude the proof, we combine Theorem  \ref{projy} with Theorem \ref{piu}. \end{proof}

\ssc{Hypothesis $(H')$ for Pseudo-Initial Times} \lab{sec63}

In this subsection, we examine the case of a pseudo-initial time.
Let us first quote from Li \cite{L} a useful result showing that the hypothesis $(H')$ is met, that is,
any $(\P,\ff)$-semimartingale is also a $(\P,\gg )$-semimartingale. Note that Theorem \ref{jacnew}
is an extension of the classic result due to Jacod \cite{JJ}, who dealt with the case of an initial time (i.e., the
density hypothesis) and the initial enlargement $\gg^*$.

\bt \lab{jacnew}
If $\tau$ is a pseudo-initial time then the hypothesis $(H')$ is satisfied by $\ff$ and
the initial enlargement $\gg^*$ (and thus also the progressive enlargement $\gg$).
\et

\begin{proof}
The arguments used in the demonstration of the theorem combine the idea of the proof of Jacod's theorem under
the standard density hypothesis (see \cite{JJ}) with the time change.
For details, the interested reader is referred to the proof of Theorem 7.2.1 in Li \cite{L}.
\end{proof}

We will now derive the $\gg$-semimartingale decomposition for a progressive enlargement with a pseudo-initial time.
Note that the $\gg$-semimartingale decomposition established in the literature under the density hypothesis
by Jeanblanc and Le Cam \cite{JL} (see also Kchia et al. \cite{KLP} who employed their Theorem \ref{kchia})
can be obtained as a special case of equality (\ref{c455}) by postulating that the increasing process $D$
in Definition \ref{densityr} is non-random. Our proof of decomposition \eqref{c455} is based on Theorem \ref{piu}.

\bt \label{piu7}
Let $\tau $ be a pseudo-initial time. If $\MM$ is a $(\P,\ff)$-local martingale
then the process $\whMM$ is a $(\P,\gg )$-local martingale where
\be \lab{c455}
\whMM_t = U_t -  \int_{(0,t \wedge \tau ]} (G_{u-})^{-1} \, d\wt C_u
+ \int_{(t \wedge s ,t]} (\mm_{s,u-})^{-1} \, d\langle U,\mm_{s,\cdot}\rangle_u \,\big|_{s= \tau}.
\ee
\et

\begin{proof}
We maintain the notation introduced in the proof of Theorem \ref{piu}.
In view of Proposition \ref{piu5} and the arguments used in the proof of Theorem \ref{piu}, it suffices to show
that, for any fixed $u \geq 0$, the process $(\mm_{u,t}\wh U^0_{u,t})_{t\geq u}$ is a $(\P,\ff)$-local martingale, where we set
(see \eqref{nexw})
\bde
\wh U^0_{u,t} =  U_{t} - Z_{u,t} - \VV_u
\ede
where in turn $\VV_u = \wt \BB_{u} + \Delta U_u + \wh U_{u,u}$, with the process $\wt \BB$ defined by
\eqref{piu3} and $Z_{u,t}$ given by the following expression  (see \eqref{nexw1})
\begin{align*}
Z_{u,t}  =  \int_{(u,t]} (\mm_{u,v-})^{-1} \, d\langle \MM,\mm_{u,\cdot}\rangle_v .
\end{align*}
By applying the integration by parts formula, we  obtain
\begin{align*}
d( \mm_{u,t} \wh U^0_{u,t} ) & =  U_{t-} \, d\mm_{s,t} + \mm_{s,t-}\, dU_{t} +
d\langle \MM,\mm_{s,\cdot}\rangle_t - Z_{u,t-} \, d\mm_{s,t} - \mm_{u,t}\, dZ_{u,t} - \VV_u \, dm_{u,t}\\
& =  U_{t-} \, d\mm_{s,t} + \mm_{s,t-}\, dU_{t} - Z_{u,t-} \, d\mm_{u,t} - \VV_u \, dm_{u,t},
\end{align*}
which is clearly a $(\P,\ff)$-local martingale for $t\geq u$.
\end{proof}

Let us assume, in addition, that there exists a positive $(\P,\ff)$-martingale
$L$ and a positive, $\ff$-adapted process $a$ such that the equality $\mm_{s,t} = L_t a_s$ holds for all $s\leq t$. Then
\bde
F_{u,t} = \int_{[0,u]} \mm_{s,t} \,dD_s = L_t\int_{[0,u]} a_s \,dD_s = K_u L_t
\ede
where $K_u = \int_{[0,u]} a_s \,dD_s$ and thus the $(\P,\ff)$-conditional distribution
of $\tau $ is completely separable. Then formula \eqref{c455} can also be deduced
from Corollary \ref{piu4}.

\bcor \label{munn2}
Let $\tau $ be a pseudo-initial time. Assume that there exists a positive $(\P,\ff)$-martingale
$L$ and a positive, $\ff$-adapted process $a$ such that  $\mm_{s,t} = L_t a_s$ for all $s\leq t$.
Then the $(\P,\ff)$-conditional distribution $F_{u,t}$ is completely separable. Moreover, if
$\MM$ is a $(\P,\ff)$-local martingale then the process $\whMM$ given by \eqref{c455}
is a $(\P,\gg )$-local martingale.
\ecor

\begin{proof}
It suffices to apply Corollary \ref{piu4} and observe that, for any fixed $s$, the equality
$\mm_{s,u}= a_s L_u$ holds. Hence formulae  \eqref{iiuu} and \eqref{c455} are equivalent.
\end{proof}

\section{Applications to Financial Mathematics}  \label{sec7}

In this final section, we present two applications of some of our general results established in preceding
sections to particular problems arising in the context of financial modeling.

\subsection{Arbitrage Free Markets Models}   \label{sec7.1}

In the recent paper by Coculescu et al. \cite{CJN}, the existence of an equivalent probability measure
under which the immersion property holds was shown to be a sufficient condition for a market model with enlarged
filtration to be arbitrage-free, provided that the underlying market model based on the filtration $\ff$ enjoyed this
property. In Proposition \ref{cfre}, we will show that if $\tau $ satisfies the hypothesis $(H\!P)$ (or, more precisely,
the complete separability of the conditional distribution $F_{u,t}$ holds) then, under mild
technical assumptions, the result from \cite{CJN} can be applied to the progressive enlargement~$\gg$.

\sssc{Immersion Property under an Equivalent Probability Measure} \lab{coc1}

Before studying the case of a progressive enlargement, we will first summarize briefly some results from Coculescu et al.~\cite{CJN} and make pertinent comments.
Suppose that a probability space $(\Omega, \F , \P)$ is endowed with arbitrary filtrations $\ff$ and $\kk$ such that $\ff \subset \kk$.
Let $\IPP$ stand for the class of all probability measures $\Q$ equivalent to $\P$ on $(\Omega , \F)$ such that
$\ff$ and $\kk$ satisfy the immersion property under $\Q$. Recall that the immersion property for $\ff $ and $\kk$ under $\Q$
stipulates that any $(\Q,\ff )$-local martingale is a $(\Q,\kk )$-local martingale.
The following lemma was established in \cite{CJN}.

\bl  \lab{lmc1x}
Assume that  $\QT \in \IPP$ and $\Q $ is a probability measure equivalent to $\QT$
on $(\Omega, \F )$ such that $\frac{d\Q }{d\QT }$ is $\F_{\infty}$-measurable. Then $\Q$ belongs to $\IPP$.
\el

\begin{proof}
Let $M$ be an $(\ff, \Q)$-martingale. We wish to show that
$M$ is also a $(\Q,\gg)$-martingale. This is equivalent to the property that
$M \eta $ is a  $(\QT ,\gg )$-martingale where $\eta_t := \frac{d\Q }{d\QT }|_{\G_t}$.
Since $M$ and $\eta$ are $\ff$-adapted processes and  $\QT \in \IPP$, it suffices to
show that $M \eta $ is a $(\QT , \ff )$-martingale.
To this end, we observe that $\eta$ is an $\ff$-adapted $(\QT ,\gg)$-martingale and thus an $(\QT, \ff)$-martingale.
Hence the Bayes formula yields, for all $0 \leq t <s$,
$$
\EQT (M_s \eta_s \,|\, \F_t) = \EQ (M_s \,|\, \F_t)\,\EQT ( \eta_s \,|\, \F_t) = M_t \eta_t ,
$$
so that $M$ is a $(\Q,\gg)$-martingale. By the usual localization argument, the proof can be extended to local martingales.
\end{proof}

\bl  \lab{lmc2x}
Assume that the class $\IPP$ is non-empty. Then for every $\QT \in \IPP$ there exists a probability measure
$\Q$ equivalent to $\QT$ on $(\Omega, \F )$ and such that the following conditions are met: \hfill \break
(i) the Radon-Nikod\'ym density process $\eta_t := \frac{d\Q }{d\QT }|_{\G_t}$ is $\ff$-adapted, \hfill \break
(ii) the probability measures $\Q$ and $\P $ coincide on $\ff $,  \hfill \break
(iii) $\Q $ belongs to $\IPP$, \hfill \break
(iv) every $(\P,\ff)$-local martingale is a $(\Q,\kk)$-local martingale.
\el

\begin{proof}
{(i)} Let $\rho_{\infty} = d\P /d\QT $ and let $\eta_{\infty} = \EQT (\rho_{\infty} \,|\, \F_{\infty})$.
We define the probability measure $\Q$ equivalent to $\QT$ on $(\Omega, \G )$ by setting
\be \lab{cder}
\frac{d\Q }{d\QT } = \eta_{\infty}.
\ee
Then
\bde
\eta_t := \frac{d\Q }{d\QT }\Big|_{\G_t}= \EQT ( \eta_{\infty} \,|\, \G_t )
= \EQT ( \eta_{\infty} \,|\, \F_t ) =  \EQT ( \rho_{\infty} \,|\, \F_t ) = \frac{d\Q }{d\QT }\Big|_{\F_t}
\ede
where the third equality holds since $\QT \in \IPP$ and $\eta_{\infty}$ is $\F_{\infty}$-measurable.
We thus see that (i) holds.

\noindent {(ii)} Furthermore,
$$
\frac{d\Q }{d\P}\Big|_{\F_t} 
= \frac{ \frac{d\Q }{d\QT }|_{\F_t} }{ \frac{d\P }{d\QT }|_{\F_t} }
= \frac{\EQT ( \rho_{\infty} \,|\, \F_t )}{\EQT ( \rho_{\infty} \,|\, \F_t )} =1
$$
and thus (ii) is satisfied.

\noindent {(iii)} The third statement follows from (\ref{cder}) and Lemma \ref{lmc1x}.
It is also worth noting that
$$
\frac{d\Q }{d\P}\Big|_{\G_t} =
 \frac{ \frac{d\Q }{d\QT }|_{\G_t} }{ \frac{d\P }{d\QT }|_{\G_t} }
=  \frac{\EQT ( \rho_{\infty} \,|\, \F_t )}{\EQT ( \rho_{\infty} \,|\, \G_t )},
$$
where the right-hand side can be checked to be a $(\P,\kk)$-martingale.

\noindent {(iv)} Since $\Q = \P $ on $\ff$, any $(\P,\ff)$-local martingale is a $(\Q,\ff)$-local martingale.
From part (iii), we know that $\Q \in \IPP$ and thus
any $(\Q,\ff )$-martingale is also a $(\Q,\kk)$-martingale.
Hence any $(\P,\ff)$-local martingale is also a $(\Q,\kk)$-martingale.
\end{proof}

Note that the probability measure $\P$ in Lemma \ref{lmc2x} can be replaced by any probability measure
equivalent to $\P$. This means that the assumption that the class $\IPP$ is non-empty is fairly strong;
it implies that for any probability measure $\P'$ there exists a probability measure $\Q'$
equivalent to $\P'$ on $(\Omega , \F)$, coinciding
with $\P'$ on $\ff$, and such that the immersion property between $\ff$ and $\kk$
holds under $\Q'$.  It is thus natural to ask under which (non-trivial) circumstances the class $\IPP$ is non-empty.
A partial answer to this question for the progressive enlargement is provided in Proposition~\ref{cfre}.

\brem \lab{honpr}
The property that the class $\IPP$ is non-empty is not satisfied by $\ff$ and the progressive
enlargement $\gg$ when $\tau $ is an $\F_{\infty}$-measurable random time (e.g., an honest time with respect to $\ff$),
unless it is an $\ff$-stopping time (so that $\ff = \gg$ and the immersion property is trivially satisfied).
\erem

The next result, which is also borrowed from \cite{CJN}, provides a complete characterization of non-emptiness of the class $\IPP$.

\bp \lab{pc2x}
The following conditions are equivalent: \hfill \break
(i) the class $\IPP$ is non-empty, \hfill \break
(ii) there exists a probability measure $\Q$ equivalent to $\P$ on $(\Omega , \F)$
such that every $(\P,\ff)$-local martingale is a $(\Q,\kk)$-local martingale.
\ep

\begin{proof}
\noindent (i) $\Rightarrow$ (ii) This implication is an immediate consequence of part (iv) in Lemma  \ref{lmc2x}.

\vskip 5 pt
\noindent (ii) $\Rightarrow$ (i)  Assume that there exists a probability measure $\Q $ equivalent to
$\P$ on $(\Omega , \G)$ such that any $(\ff, \P)$-local martingale is a $(\Q,\kk)$-local martingale.
We claim that the hypothesis $(H)$ holds under $\Q$.
To this end, it suffices to show that $\Q = \P$ on $\ff$. Then any $(\Q ,\ff )$-local martingale is obviously
an $(\P,\ff )$-local martingale and thus, by assumption, it is also a $(\Q,\kk)$-local martingale.
This means that the hypothesis $(H)$ holds under $\Q $.

It remains to show that (ii) implies that $\Q = \P$ on $\ff$. For this purpose,
 we note that the process $\eta_t := \frac{d\Q }{d\P }|_{\F_t}$ is an $(\P,\ff)$-martingale.
Hence, by assumption, it is also a $(\Q,\kk )$-martingale.
Since $\eta $ is an $\ff$-adapted process, we see that it is a $(\Q ,\ff)$-martingale.
By the usual argument (see Lemma \ref{lmc1x}), this in turn is equivalent
to the property that $\eta^2$ is an $(\P,\ff)$-local martingale. Since $\eta$ and $\eta^2$ are $(\P,\ff)$-martingales,
we conclude that $\eta $ is constant process: $\eta =1$, and thus the probability measures
$\Q $ and $\P$ coincide on $\ff$.
\end{proof}

\sssc{Martingale Measures via the Hypothesis $(H)$} \lab{coc2}

Suppose now we are given a $(\P,\ff)$-semimartingale $X$ defined on a probability space $(\Omega, \F , \P)$.
Let ${\cal M}(\P,\ff)$ stand for the class of all {\it $\ff$-local martingale measures} for $X$, meaning that
a probability measure $\Q$ belongs to ${\cal M}(\P,\ff)$ whenever (i) $\Q$ equivalent to $\P$ on $(\Omega ,\F)$
and (ii) $X$ is an $(\Q,\ff)$-local martingale.
Let $\kk$ be any enlargement of the filtration $\ff$. We denote by ${\cal M} (\P,\kk)$
the class of {\it $\kk$-local martingale measures} for $X$.
In \cite{CJN}, the authors assumed that the class  ${\cal M}(\P,\ff)$ is non-empty and they searched for
sufficient conditions ensuring that the class  ${\cal M}(\P,\kk)$ is non-empty as well.

One possibility is to postulate that the immersion property holds under some
probability measure $\Q$ equivalent to $\P$ and to infer that it
is also valid under some $\ff$-local martingale measure. Obviously, any  $\ff$-local martingale measure under
which the immersion property holds is also a $\kk$-local martingale measure.
A particular example of an $\ff$-local martingale measure under
which the the immersion property holds can be produced using Lemma  \ref{lmc2x}, leading to
the following result, also due to Coculescu et al. \cite{CJN} (see Corollary 4.6 therein).

\bp \lab{pc1x}
(i) The classes ${\cal M}(\P,\ff)$ and $\IPP$ are non-empty if and only if the set ${\cal M} (\P,\ff) \cap \IPP$ is non-empty. \hfill \break
(ii) If the classes ${\cal M}(\P,\ff)$ and $\IPP$ are non-empty then the class ${\cal M} (\P,\kk)$ is non-empty.
\ep

\begin{proof}
(i) Assume that ${\cal M} (\P,\ff)$ and $\IPP$ are non-empty. Let $\QT \in \IPP$ and let  $\P' \in {\cal M} (\P,\ff)$.
Let $\Q$ be defined as in Lemma  \ref{lmc2x} with $\P$ replaced by $\P'$ (of course, $\IPP=\IPPS $ since $\P$ is equivalent to $\P'$).
Then $\Q \in \IPP$ by part (iii) in Lemma \ref{lmc2x}. Also, by part (ii) in Lemma \ref{lmc2x},
the probability measures $\Q$ and $\P'$ coincide on $\ff$ and thus $\Q \in {\cal M} (\P,\ff)$. Hence $\Q \in {\cal M} (\P,\ff) \cap \IPP$
and thus the class ${\cal M} (\P,\ff) \cap \IPP$ is non-empty. The converse implication is trivial.

\noindent (ii) In view of part (i), it suffices to show that any probability
measure $\Q$ belonging to ${\cal M} (\P,\ff) \cap \IPP$ is in ${\cal M} (\kk ,\P)$.
Indeed, any $(\Q,\ff )$-local martingale is a $(\Q,\kk )$-local martingale (since the immersion property holds under $\Q$)
and $X$ is an $(\Q,\ff )$-local martingale (since $\Q$ is an $\ff$-local martingale measure),
so that $X$ is a $(\Q,\kk )$-local martingale, as required.
\end{proof}

\brem An inspection of the proof of Proposition \ref{pc1x} (i.e. Corollary 4.6 in \cite{CJN}) shows that the process $X$ plays no essential role
(except, of course, for the assumption that the class  ${\cal M}(\P,\ff)$ is non-empty).
Note also that the probability measure $\Q$ constructed in this proof
has the property that every $(\P',\ff)$-local martingale is also an
$(\Q,\ff )$-local martingale (since $\P'$ and $\Q$ coincide on $\ff$)
and thus a $(\Q,\kk )$-local martingale (since $\Q \in \IPPS$).
\erem

\brem
The class ${\cal M}(\P,\ff)$ can be replaced in part (i) of Proposition \ref{pc1x}
by any subset ${\cal P}$ of probability measures equivalent to $\P$ and such that:
if $\P' \in {\cal P}$ and $\P'' = \P' $ on $\ff $ then $\P'' \in {\cal P}$.
\erem

To summarize the conclusions from  Coculescu et al. \cite{CJN}, the following conditions are equivalent: \hfill \break
(1)  ${\cal M} (\P,\ff)$ and $\IPP$ are non-empty, \hfill \break
(2)  ${\cal M} (\P,\ff) \cap \IPP$ is non-empty, \hfill \break
(3)  ${\cal M} (\P,\ff)$ is non-empty and there exists a probability measure $\Q$ equivalent to $\P$ on $(\Omega , \F)$
such that every $(\P,\ff)$-local martingale is a $(\Q,\kk )$-local martingale.

In view of Proposition \ref{pc1x}, any of conditions (1)--(3) implies that the class ${\cal M} (\P,\kk )$ is non-empty.
It is thus natural to refer to any of conditions (1)--(3) is the {\it no-arbitrage condition} for the market model with
an enlarged filtration $\kk$.

\sssc{Martingale Measures for the Progressive Enlargement} \lab{coc2x}

We now consider the case where $\kk = \gg$ is the progressive enlargement of $\ff$.
Suppose that the hypothesis $(H)$ is not satisfied by the $(\P,\ff)$-conditional distribution $F_{u,t}$ of a random time $\tau $.
It is then natural to ask whether there exists a probability measure $\bar \P$, which is equivalent to $\P$ on $(\Omega , \F)$ and such that
the $(\bar \P,\ff)$-conditional distribution $\bar F_{u,t}$ of $\tau $ satisfies the hypothesis $(H)$.
Equivalently, we ask whether there exists a probability measure $\bar \P$ equivalent to $\P$ and such that, for all $0 \leq u \leq t$,
\be \lab{xdee}
\bar F_{u,u} := \cP{\bar \P}{\tau \leq u}{\F_u} = \cP{\bar \P}{\tau \leq u}{\F_t} =: \bar F_{u,t}.
\ee
Under the density hypothesis, the answer to this question is known to be positive (see
El Karoui et al. \cite{EJJ} and Grorud and Pontier \cite{GP}). By contrast, when $\tau $ is assumed to be an honest time
then this property never holds, unless $\tau $ is an $\ff$-stopping time (see Remark \ref{honpr}).

In this subsection, we work under the standing assumption that the $(\P,\ff)$-conditional distribution of $\tau$ is separable
and $F_0=0$, so that $\tau $ is a pseudo-honest time (see Remark \ref{uiui}). Note, however, that the case of the classic honest time is not covered by foregoing
results, since we also assume from now on that $F_{u,t}>0$ for all $0 < u \leq t$. Recall also that the complete separability implies
that $\tau $ is a pseudo-honest time.

\bl
Assume that $(\P,\ff)$-conditional distribution of $\tau $ is completely separable, so that
$F_{u,t} = K_u L_t$ for $0 \leq u \leq t$. Let the process $(Z^\gg_t)_{t\geq 0}$ be given by
\be \label{HPtoH}
Z^\gg_t = \wt Z_t\I_\seq{\tau > t} + \wh Z_{\tau,t}\I_\seq{\tau \leq t}
\ee
where $\wh Z_{u,t} = \frac{F_{u,u}}{F_{u,t}}$ and
\bde
\wt Z_t= (G_t)^{-1} \bigg(1- \int_{(0,t]} \wh Z_{u,t} \, dF_{u,t} \bigg)
= (G_t)^{-1}\bigg(1- \cE{\P}{\wh Z_{\tau,t}\I_\seq{\tau \leq t}}{\F_t} \bigg).
\ede
Then the process $Z^\gg$ is a $(\P,\gg)$-local martingale.
\el

\begin{proof}
The proof relies on an application of Corollary \ref{xeed} and Theorem \ref{piu1}. First, $\wh Z_{u,u} = 1$ and thus it is trivially an $\ff$-predictable process.
Next, we need to check condition (i) in Corollary \ref{xeed},
which now reads: the process $(W_{t})_{t \geq 0}$ is a $(\P,\ff)$-local martingale where
\be
W_t = \wt Z_t G_t + \int_{(0,t]} \wh Z_{u,u} \, dF_u  = \wt Z_t G_t + F_t.
\ee
We observe that
\begin{align*}
\wt Z_t G_t + F_t &= 1 + F_t - \int_{(0,t]} \wh Z_{u,t} \, dF_{u,t} \\
&= 1+ L_tK_t - \int_{(0,t]} L_u \, dK_u = 1 + L_0K_0 + \int_{(0,t]} K_{u-} \, dL_u ,
\end{align*}
so that the process $W$ is indeed a $(\P,\ff)$-local martingale.
Finally, we need to check condition (ii) in Corollary \ref{xeed}, which takes here the following form:
for every $u>0$, the process $(W_{u,t} = L_t (\wh Z_{u,t} - \wh Z_{u,u}))_{t \geq u}$
is a $(\P,\ff)$-local martingale. To this end, we note that
\be
W_{u,t} = L_t\frac{F_{u,u}}{F_{u,t}}-L_t \wh Z_{u,u} = L_t \frac{K_uL_u}{K_uL_t}-L_t \wh Z_{u,u} = L_u -
L_t \wh Z_{u,u} , \label{c2}
\ee
which is, obviously, a $(\P,\ff)$-local martingale for $t\geq u$.  In view of Theorem \ref{piu1}, we conclude
that $Z^\gg$ is a $(\P,\gg)$-local martingale.
\end{proof}

We will also need the following simple lemma.

\bl \lab{tgtg}
Assume that $\tau $ is a pseudo-honest time. Then for any $\ff$-adapted, $\P$-integrable process $X$ we have that,
for every $s \leq t$,
\be \lab{wsz}
F_{s,t} \, \EP ( X_{\tau } \I_{\{ \tau \leq s\}} \,|\, \F_s ) = F_{s,s} \, \EP ( X_{\tau } \I_{\{ \tau \leq s\}} \,|\, \F_t ).
\ee
\el

\begin{proof}
Note that for $X=1$ equality (\ref{wsz}) is trivially satisfied.  In general, it suffices to
consider an elementary $\ff$-adapted process of the form $X_t = \I_A \I_{[u,\infty )}(t)$ for a fixed, but arbitrary, $u \geq 0$ and
any event $A \in \F_u$. Obviously, both sides of (\ref{wsz}) vanish when $u >s$. For any $u \leq s$, we obtain
\begin{align*}
&F_{s,t} \, \EP ( \I_A \I_{[u,\infty )}(\tau ) \I_{\{ \tau \leq s\}} \,|\, \F_s ) =\I_A F_{s,t} \, \EP (  \I_{\{ \tau \leq s\}}- \I_{\{ \tau < u \}} \,|\, \F_s )
= \I_A F_{s,t}  (F_{s,s} - F_{u-,s}) \\  & =  \I_A ( F_{s,t}F_{s,s} - F_{u-,t}F_{s,s})
= \I_A F_{s,s} ( F_{s,t} - F_{u-,t}) = \I_A F_{s,s} \, \EP (  \I_{\{ \tau \leq s\}}- \I_{\{ \tau < u \}} \,|\, \F_t )
\\ &= F_{s,s} \, \EP ( \I_A \I_{[u,\infty )}(\tau ) \I_{\{ \tau \leq s\}} \,|\, \F_t )
\end{align*}
where we used condition (\ref{hpx}) in the third equality.
\end{proof}

\bp \lab{cfre}
Assume that: \hfill \break
(i) the $(\P,\ff)$-conditional distribution of a random time $\tau$ is completely separable and $F_0=0$, \hfill \break
(ii) the process $Z^\gg$ given by formula (\ref{HPtoH}) is a positive $(\P, \gg)$-martingale
such that $\EP (Z^\gg_t\,|\,\F_t) =1$ for every $t \in \rr_+$. \hfill \break
Then there exists an equivalent probability measure $\bar \P $ such  $\bar \P = \P$ on $\ff$ and the hypothesis
 $(H)$ holds under $\bar \P$.
\ep

\begin{proof}
It suffices to show that (\ref{xdee}) holds under $\bar \P$, where the probability measure $\bar \P$
is defined on $\G_t$ by $d\bar \P|_{\G_t} = Z^\gg_t \, d\P|_{\G_t}$. We observe that, for all $0 \leq u \leq t$,
\begin{align*}
\bar F_{u,u}  &= \cP{\P}{Z^{\gg}_u \I_\seq{\tau \leq u}}{\F_u}
= \cP{\P}{\wh Z_{\tau,u} \I_\seq{\tau \leq u}}{\F_u}
= (X_{u})^{-1}\, \cP{\P}{ X_{\tau} \I_\seq{\tau \leq u}}{\F_u} \\
&=(X_{t})^{-1}\, \cP{\P}{ X_{\tau} \I_\seq{\tau \leq u}}{\F_t} = \cP{\P}{\wh Z_{\tau,t} \I_\seq{\tau \leq u}}{\F_t}
= \cP{\P}{Z^{\gg}_t \I_\seq{\tau \leq u}}{\F_t} = \bar F_{u,t}
\end{align*}
where the fourth equality is an immediate consequence of Lemma \ref{tgtg}.
\end{proof}

An important conclusion from Proposition \ref{cfre} is that if $\tau $ is a strictly positive random time such that
the complete separability of the $(\P,\ff)$-conditional distribution holds then, under technical conditions of
Proposition \ref{cfre}, we have that $\IPP \neq \emptyset$. Therefore, part (ii) in Proposition \ref{pc1x} can be
applied in these circumstances to the progressive enlargement $\gg$.

\subsection{Information Drift}   \label{sec7.2}

For simplicity of presentation, we assume in this subsection that all $(\P,\ff)$-local martingales are continuous,
that is, assumption $({\bf C})$ is valid.
Our aim is to apply the semimartingale decompositions developed in Section \ref{sec6}
to utility maximization and information theory associated with continuous time models of financial markets,
as studied, in particular, by Ankirchner and Imkeller~\cite{AI}.
We will need the following definition borrowed from Ankirchner and Imkeller~\cite{AI} (see Definition 1.1 in \cite{AI}).
 In what follows, the process $X$ can be interpreted as
the discounted price of a risky asset. The interested reader is
referred to \cite{AI} for the compelling rationale for this definition in the context of maximization
of the expected logarithmic utility from a portfolio's wealth when trading in $X$.

\bd  \lab{deffin}
A filtration $\ff$ is said to be a {\it finite utility filtration} for a process $X$ whenever $X$ is a $(\P,\ff)$-semimartingale
with the semimartingale decomposition of the form
\be\label{v}
X_t = U_t + \int_{(0,t]} \phi_u \, d\langle U,U \rangle_u
\ee
for some $(\P,\ff)$-local martingale $U$ and an $\ff$-predictable process $\phi$.
\ed

\brem
It was shown by Delbaen and Schachermayer \cite{DeS} that for a locally bounded semimartingale $X$,
the existence of an equivalent local martingale measure for $X$ (or, equivalently, the property that
$X$ satisfies the No Free Lunch with Vanishing Risk (NFLVR) condition) implies that
the decomposition of the process $X$ must be of the form given in equation \eqref{v}.
\erem

From now on, we work on a filtered probability space $(\Omega, \F, \ff , \P)$ with a filtration $\ff$
satisfying the usual conditions and such that $\ff$ is a finite utility filtration for a given process $X$
so that decomposition \eqref{v} holds for some $(\P,\ff)$-martingale $U$.
The following definition comes from Imkeller \cite{I} (see also Definition 1.2 in \cite{AI}).

\bd \label{drift}
Let $\kk$ be any enlargement of the filtration $\ff $. The $\kk$-predictable process $\psi$ such that the process
\bde
U_t-\int_{(0,t]} \psi_u \, d\langle U,U \rangle_u
\ede
is a $(\P,\kk)$-local martingale is called the {\it information drift} of $\kk$ with respect to $\ff$.
\ed

\brem
It was shown in \cite{AI} (see Proposition 1.2 therein) that the difference
of the maximal expected logarithmic utilities when trading in $X$ is based on
two different filtrations $\ff \subset \kk$ depends only on the
information drift of $\kk$ with respect to $\ff$.
\erem

Under the standing assumption $({\bf C})$  all $(\P,\ff)$-local martingales encountered in what follows are locally
bounded, and thus locally square-integrable, so that we are in a position to apply the Kunita-Watanabe decomposition theorem
(see, for instance, Protter \cite{PP}).
We aim to show that the progressive enlargement of $\ff$ with a pseudo-initial time  $\tau$
is once again a finite utility filtration for $X$. We will also compute the information drift of
the progressive enlargement $\gg$ with respect to $\ff$.

\bp
If $\ff$ is a finite utility filtration for $X$ and $\tau$ is a pseudo-initial time then
the progressive enlargement $\gg$ is a finite utility filtration for $X$. Furthermore,
the information drift of $\gg$ with respect to $\ff$ is given by the following expression
\be \label{gmmac}
\psi_u = \I_\seq{\tau \geq u }(G_u)^{-1}\eta_u
 + \I_{\{ \tau < u \}}  (m_{\tau,u})^{-1} \xi_{\tau,u}
\ee
where the $\ff$-predictable processes $(\xi_{s,u})_{s \leq u}$ are defined by the Kunita-Watanabe decompositions
\be\label{kw2}
m_{s,t} = \int_{(s,t]} \xi_{s,u}\, dU_u + L^s_t
\ee
where $L^s$ is a family of $(\P,\ff)$-local martingales strongly orthogonal to $U$ for all $u$.
\ep

\begin{proof}
By the standing assumption, the filtration $\ff$ is a finite utility filtration for the process $X$. Therefore,
the process $X$ admits the following  $(\P,\ff)$-semimartingale decomposition (see \eqref{v})
\bde
X_t = U_t + \int_{(0,t]} \phi_u \, d\langle U,U \rangle_u
\ede
where $U$ is a $(\P,\ff)$-local martingale and $\phi $ is an $\ff$-predictable process. To establish the first assertion,
it suffices to show that the process $X$ is a $(\P,\gg)$-semimartingale with the following decomposition
\be\label{sg}
X_t = \whMM_t + \int_{(0,t]} \phi^*_u \,d\langle \whMM,\whMM \rangle_u
\ee
where $\whMM$ is some $(\P,\gg)$-martingale and $\phi^*$ is some $\gg$-predictable process. Using Corollary \ref{munn2}
and the assumption that all $(\P,\ff)$-local martingales are continuous (so that $\wt C = \langle U , M \rangle)$, we deduce
that the process $\whMM$, which is given by the expression
\be\label{info1}
\whMM_t = U_t - \int_{(0,t]}  \I_\seq{\tau \geq u }(G_u)^{-1} \, d\langle U,M \rangle_u
 -  \int_{(t \wedge s ,t]} (\mm_{s,u-})^{-1} \, d\langle U,\mm_{s,\cdot}\rangle_u \,\big|_{s= \tau},
\ee
is a $(\P,\gg)$-local martingale. Recall that we denote by $M$ the $(\P,\ff)$-local martingale
 appearing in the Doob-Meyer decomposition of $G$.
An application of the Kunita-Watanabe decomposition theorem to $M$ and $U$, after suitable localization
if required, gives
\be\label{kw1}
M_t = \int_{(0,t]} \eta_u \,dU_u + \wh L_t
\ee
where $\eta$ is some $\ff$-predictable process and a square-integrable $(\P,\ff)$-martingale $\wh L$ is strongly orthogonal to $U$.
It thus follows immediately from \eqref{kw1} that $d\langle U,M\rangle_u = \eta_u \, d\langle U,U\rangle_u$.
In the next step, we focus on the $(\P,\ff)$-martingale $(m_{s,t})_{t\geq s}$, for any fixed $s \in \rr_+$.
Using once again the Kunita-Watanabe decomposition theorem, we deduce that \eqref{kw2} holds for all $s \leq t$
where $(\xi_{s,u})_{u \geq s}$ is a family of $\ff$-predictable processes parametrized by $s$ and $L^s$ a
family of $(\P,\ff)$-local martingales strongly orthogonal to $U$ for every $s$.
Consequently, by combining  \eqref{kw2}, \eqref{info1} and \eqref{kw1}, we arrive at the following equalities
\begin{align*}
\whMM_t & =  U_t - \int_{(0,t]}  \I_\seq{\tau \geq u }(G_u)^{-1}\eta_u \, d\langle U,U\rangle_u
 + \int_{(0,t]} \I_{\{ \tau < u \}}  (m_{\tau,u})^{-1} \xi_{\tau,u}\, d\langle U,U\rangle_u\\
	 & = U_t - \int_{(0,t]}  \Big( \I_\seq{\tau \geq u }(G_u)^{-1}\eta_u
 + \I_{\{ \tau < u \}}  (m_{\tau,u})^{-1} \xi_{\tau,u} \Big) \, d\langle U,U\rangle_u .
\end{align*}
It is clear that the equality $\langle U,U\rangle =\langle \whMM,\whMM\rangle$ holds.
Therefore, the canonical $(\P,\gg)$-semimartingale decomposition of $U$ reads
\be\label{ug}
U_t = \whMM_t + \int_{(0,t]}  \psi_u \, d\langle \whMM,\whMM\rangle_u
\ee
where the $\gg$-predictable process $\psi$ is given by the following expression
\bde
\psi_u = \I_\seq{\tau \geq u }(G_u)^{-1}\eta_u
 + \I_{\{ \tau < u \}}  (m_{\tau,u})^{-1} \xi_{\tau,u}.
\ede
It is now easy to see that the $(\P,\gg)$-semimartingale decomposition of $X$ has indeed the desired form
 \eqref{sg} with $\phi^* = \phi + \psi$.
To derive equality (\ref{gmmac}), it is enough to employ Definition \ref{drift}. The asserted
formula follows directly from representation \eqref{ug}
and the fact that $\langle \whMM,\whMM \rangle = \langle U,U\rangle$.
\end{proof}



\vskip 15 pt

\section{Appendix: An Overview of Classic Results}

The goal of the appendix is to provide a succinct overview of classic results regarding the
semimartingale decomposition for the progressive enlargement of filtration with the special
emphasis on the case of an honest time. In particular, we give here the proofs of Theorems \ref{projy}
and \ref{thhon}. For more information and further results, the interested reader is referred to Jeulin \cite{J1,J2},
Jeulin and Yor \cite{JY1,JY2,JY3} and Yor \cite{Y1,Y2}.

Let $\tau $ be any strictly positive random time on a filtered probability space $(\Omega , {\cal F} , \ff , \P)$. We use the following notation
\bde
\II = \I_{\llb 0,\tau \rrb}(t), \quad \JJ = \I_{\rrb \tau , \infty \llb}(t),
\ede
\bde
\tII =  \I_{\llb 0,\tau \llb}(t), \quad \tJJ = \I_{\llb \tau , \infty \llb}(t).
\ede
We use the standard notation for the (dual) $\ff$-optional and $\ff$-predictable and  projections.
In particular, the dual $\ff$-optional projection $\tJJ^o$ (resp. the dual $\ff$-predictable projection $\tJJ^p$) is a positive, increasing, $\ff$-optional (resp. $\ff$-predictable) process. We denote $(X \cdot Y)_t = \int_{(0,t] } X_u \, dY_u.$

\brem
Our usual notation reads: $H = \tJJ $. For optional projections, we denote
\bde
G  = {}^o(1-H) = {}^o \tII , \quad F = {}^o H = {}^o \tJJ = 1- G .
\ede
It is known that ${}^p \II = G_-$  (see Jeulin \cite{J1}). Hence ${}^p \JJ = {}^p(1- \II )= 1-G_- = F_-$. Note that ${}^p F \ne F_-$, in general,
since it may happen that ${}^p \tJJ \ne {}^p \JJ$ or, equivalently, that ${}^p\I_{\llb \tau \rrb} \ne 0$.
\erem

\bl \lab{rr5}
Let $N$ be a bounded $(\P,\ff)$-martingale with $N_0=0$. Then
\begin{align} \lab{casea}
\EP (N_{\tau }) &= \EP ( (N \cdot \tJJ)_{\infty} ) = \EP ( (N \cdot \tJJ^o)_{\infty})
\stackrel{\rm Ito}= \EP ( N_{\infty} \tJJ^o_{\infty}) \\
& = \EP ( N_{\infty} \bar M_{\infty}) =
 \EP ( [N, \bar M]_{\infty }) = \EP ( \langle N, \bar M \rangle_{\infty }) \nonumber
\end{align}
where we set $\bar M_t = \EP (  \tJJ^o_{\infty } \,|\, \F_t )$ so that  $\bar M_{\infty} = \tJJ^o_{\infty }$.
\el

\newpage

\brem   \lab{remvv1x}
One can check that ${}^o\tII  = \bar M - \tJJ^o$ and thus the equality $G = \bar M - \tJJ^o$ is
an $\ff$-optional decomposition of $G$. To this end, it suffices to show that for any
increasing, $\ff$-adapted process $B$ with $B_0=0$ we have that
\be \lab{uuu}
\EP (( ( \bar M - \tJJ^o ) \cdot B )_{\infty}) = \EP ((  {}^o\tII  \cdot B )_{\infty}).
\ee
To show that \eqref{uuu} is valid, we observe that, on the one hand,
\bde
\EP ((  {}^o\tII  \cdot B )_{\infty})= \EP (( \tII  \cdot B )_{\infty}) = \EP ( B_{\tau -}).
\ede
On the other hand, by setting  $\bar Z = \bar M - \tJJ^o$, we obtain from the It\^o formula (note that $\bar Z_{\infty}= 0$)
\bde
\EP (( \bar Z \cdot B )_{\infty}) = \EP ( (\bar ZB)_{\infty} - (B_- \cdot \bar M )_{\infty } + (B_- \cdot \tJJ^o )_{\infty })
 = \EP ( (B_- \cdot \tJJ )_{\infty }) = \EP ( B_{\tau -}).
\ede
We conclude that \eqref{uuu} is satisfied.
\erem

\brem One can show that the $(\P,\ff)$-martingale $\bar M$ belongs to the class BMO
and it is the unique BMO $(\P,\ff)$-martingale for which the equality $ \EP (N_{\tau })= \EP ( \langle N, \bar M \rangle_{\infty })$
holds for any bounded $(\P,\ff)$-martingale $N$.
\erem

\bl \lab{rr6}
Let $N$ be a bounded $(\P,\ff)$-martingale with $N_0=0$.
Then
\begin{align} \lab{caseb}
\EP (N_{\tau -}) &= \EP ((N_{-})_{\tau }) = \EP ( (N_- \cdot \tJJ)_{\infty} ) = \EP ( (N_- \cdot \tJJ^p)_{\infty})
\stackrel{\rm Ito} = \EP ( N_{\infty} \tJJ^p_{\infty}) \\
&= \EP ( N_{\infty} \wt M_{\infty}) = \EP ( [N, \wt M]_{\infty }) = \EP ( \langle N, \wt M\rangle_{\infty })  \nonumber
\end{align}
where we set $\wt M_t = \EP (  \tJJ^p_{\infty } \,|\, \F_t )$ so that $\wt M_{\infty} = \tJJ^p_{\infty }$.
\el

\brem \lab{remvvv}
We wish to check that $G = \wt M - \tJJ^p$ is the Doob-Meyer decomposition of $G$. It is enough to prove
that the equality ${}^o\tII  =\wt M - \tJJ^p$ is satisfied.  As in Remark \ref{remvv1x}, it suffices to show that for any
increasing, $\ff$-adapted process $B$ with $B_0=0$ we have that
\be \lab{xuuu}
\EP (( ( \wt M - \tJJ^p ) \cdot B )_{\infty}) = \EP ((  {}^o\tII  \cdot B )_{\infty}).
\ee
On the one hand, as before we have that
\bde
\EP ((  {}^o\tII  \cdot B )_{\infty})= \EP (( \tII  \cdot B )_{\infty}) = \EP ( B_{\tau -}).
\ede
On the other hand, by setting  $Z =\wt M - \tJJ^p$, we obtain from the It\^o formula (note that $Z_{\infty}= 0$)
\bde
\EP (( Z \cdot B )_{\infty}) = \EP ( (ZB)_{\infty} - (B_- \cdot \wt M )_{\infty } + (B_- \cdot \tJJ^p )_{\infty })
 = \EP ( (B_- \cdot \tJJ )_{\infty }) = \EP ( B_{\tau -}).
\ede
We conclude that \eqref{xuuu} is satisfied.
\erem

\brem One can show that the $(\P,\ff)$-martingale $\wt M$ belongs to the class BMO
and it is the unique BMO $(\P,\ff)$-martingale for which the equality $ \EP (N_{\tau -})= \EP ( \langle N, \wt M \rangle_{\infty })$
holds for any bounded $(\P,\ff)$-martingale $N$.
\erem

\brem
Let $G=M-A$ be the Doob-Meyer decomposition of $G$. We claim that $M = \wt M$ and $A = \tJJ^p$. To this end,
we note that $G = {}^o \tII$ and, for any random time $\tau $, the equality ${}^o \tII =\wt M-\tJJ^p$ holds,
where the $(\P,\ff)$-martingale $\wt M$ is defined by $\wt M_t := \EP (  \tJJ^p_{\infty } \,|\, \F_t )$.
Hence the uniqueness of the Doob-Meyer decomposition implies that $M = \wt M$ and $A =\tJJ^p$.
\erem

\brem
It is worth noting that we may equally well use $\tJJ^o$ and $\bar M$ in Lemma \ref{rr6}.
However, a suitable `correction term' will appear in formula \eqref{caseb}. Specifically, \eqref{caseb} will become
\bde
\EP (N_{\tau -}) =  \EP ( \langle N, \bar M\rangle_{\infty }) - \EP ( \langle N,\tJJ^o \rangle_{\infty }).
\ede
\erem

\brem
 If $\tau $ avoids $\ff$-stopping times then $\bar M = \wt M$. The avoidance property ({\bf A}) implies that ${}^o(\I_{\llb \tau \rrb}) =0$
and thus also ${}^p(\I_{\llb \tau \rrb}) = {}^p( \tJJ - J ) =0$. Hence ${}^p F ={}^p \tJJ =  {}^p \JJ = F_-$. Furthermore,
 $(\I_{\llb \tau \rrb})^o =(\I_{\llb \tau \rrb})^p =0$. More details on the avoidance property ({\bf A}) are required.
 Also, one can study the property ({\bf C}) stating that all $(\P,\ff)$-local martingales are continuous.
\erem

\ssc{Stopped Processes: Arbitrary Random Times}

Let $\kk$ be any filtration, for instance, any enlargement of $\ff$.
The following result is well known.

\bl \lab{rr7}
A bounded $\kk$-semimartingale $Y$ is a $\kk$-martingale if and only if  for any bounded $\kk$-predictable process $X$
we have
that
\bde
\EP ( (X \cdot Y)_{\infty }) = 0 .
\ede
\el

In the rest of this section, we assume that $\kk$ is any enlargement of $\ff$  such that $\tau $ is a $\kk$-stopping time and
$\kk$ is admissible prior to $\tau $. Recall that $\kk$ is {\it admissible prior to} $\tau$ if for every $t \in \rr_+$
\bde
\K_t \cap \{ \tau > t \} = \F_t \cap \{ \tau > t \},
\ede
that is: for every $t \in \rr_+$ and $A_t \in \K_t$ there exists $\wt A_t \in \F_t$ such that $A_t \cap \{ \tau > t \} = \wt A_t \cap \{ \tau > t \}$.

\bl
For any $\kk$-predictable process $X$ there exists a unique $\ff$-predictable process $K$ such that
$X \I_{\llb 0,\tau\rrb} = K \I_{\llb 0,\tau\rrb}$.
If $X$ is bounded then one make take $K$ bounded by the same constant.
\el


\bl \lab{rry4}
Let $X$ be a $\kk$-predictable process. Then
\be
{}^p (X \II ) = {}^p (K\II )  = {}^p \II K . 
\ee
\el

We assume that $Y$ is a bounded $(\P,\ff)$-martingale and
we take for granted that  the stopped process $Y^{\tau }$ is a special $(\P,\kk )$-semimartingale.
We only address the problem of explicit computation of the $(\P,\kk )$-canonical
decomposition of $Y^{\tau }$.
Propositions \ref{ppp1} and \ref{ppp2} furnish two alternative solutions to this problem,
originally due to Jeulin \cite{J2} (Proposition 4.16 in \cite{J2})
and Jeulin and Yor \cite{JY1} (Theorem 1 in \cite{JY1}), respectively.


\bp \lab{ppp1}
Assume that $Y$ is a bounded $(\P,\ff)$-martingale. The $(\P,\kk)$-canonical decomposition of the stopped process $Y^{\tau }$
reads $Y^\tau  = \widehat M + \widehat A$ where the $\kk$-predictable process of finite variation $\widehat A$ equals
\be
 \widehat A_t = \Big( \frac{\II}{\IIp} \cdot \langle Y , \bar M \rangle \Big)_t
 =\int_0^t \frac{\II_u}{\IIp_u} \, d \langle Y , \bar M \rangle_u .
\ee
\ep

\begin{proof}
The proof hinges on Lemma \ref{rr7}.
Let thus $X$ be any bounded $\kk$-predictable process. Recall that we denote by $K$ the unique bounded $\ff$-predictable process such that
$X \I_{\llb 0,\tau\rrb} = K \I_{\llb 0,\tau\rrb}$. Hence
\bde
\EP ( ( X \cdot Y^{\tau })_{\infty }) = \EP ( ( K \cdot Y^{\tau } )_{\infty })= \EP ( ( K \cdot Y )_{\tau }).
\ede
By applying Lemma \ref{rr5} to the bounded $(\P,\ff)$-martingale $N = K \cdot Y$, we obtain
\bde
\EP ( ( K \cdot Y )_{\tau }) = \EP ( \langle K \cdot Y, \bar M \rangle_{\infty }) =
    \EP ( (K \cdot \langle Y, \bar M \rangle )_{\infty })
\ede
and Lemma \ref{rry4} yields
\bde
\EP ( (K \cdot \langle Y, \bar M \rangle )_{\infty }) =
\EP \Big( \Big(K  \frac{\II}{\IIp} \cdot \langle Y, \bar M \rangle \Big)_{\infty }\Big)=
\EP \Big( \Big(X  \frac{\II}{\IIp} \cdot \langle Y, \bar M \rangle \Big)_{\infty }\Big)
= \EP ( ( X \cdot \wh A )_{\infty }).
\ede
We conclude that for any bounded $\kk$-predictable process $X$ the following equality holds
\bde
\EP ( ( X \cdot Y^{\tau })_{\infty }) = \EP ( ( X \cdot \wh A )_{\infty }).
\ede
Therefore, in view of Lemma \ref{rr7}, the process $\wh M = Y^{\tau } - \wh A$ is a $\kk$-local martingale.
It is also clear that $\wh A$ is a $\kk$-predictable process of finite variation and thus the equality
$Y^\tau  = \widehat M + \widehat A$ is the $(\P,\kk)$-canonical decomposition of the stopped process $Y^{\tau }$.
\end{proof}

\bp \lab{ppp2}
Assume that $Y$ is a bounded $(\P,\ff)$-martingale. The $(\P,\kk)$-canonical decomposition of the stopped process $Y^{\tau }$
reads $Y^\tau  = \wt M + \wt A$ where the $\kk$-predictable process of finite variation $\wt A$ equals
\be
 \widetilde A_t = \Big( \frac{\II}{\IIp} \cdot ( \langle Y ,  M \rangle + \check Y) \Big)_t
 =\int_0^t \frac{\II_u}{\IIp_u} \, d  ( \langle Y ,  M \rangle_u + \check Y_u )
\ee
where we denote $\check Y = ( \Delta Y \cdot \tJJ)^p$.
\ep

\begin{proof}
Let $X$ be any bounded $\kk$-predictable process. As before, we denote by $K$ the unique bounded $\ff$-predictable process such that
$X \I_{\llb 0,\tau\rrb} = K \I_{\llb 0,\tau \rrb }$. We note that
\bde
\EP ( ( X \cdot Y^{\tau })_{\infty }) =  \EP ( ( K \cdot Y^{\tau } )_{\infty })=
\EP ( ( K \cdot Y )_{\tau -}) + \EP ( K \cdot ( \Delta Y \cdot \tJJ )_{\infty } ) .
\ede
By applying  Lemma \ref{rr6} to the bounded $(\P,\ff)$-martingale $N = K \cdot Y$, we obtain
\bde
\EP ( ( K \cdot Y )_{\tau -}) = \EP ( \langle K \cdot Y, M \rangle_{\infty }) =
    \EP ( (K \cdot \langle Y, M \rangle )_{\infty })
\ede
and Lemma \ref{rry4} gives
\bde
\EP ( (K \cdot \langle Y, M \rangle )_{\infty }) =
\EP \Big( \Big(K  \frac{\II}{\IIp} \cdot \langle Y, M \rangle \Big)_{\infty }\Big)=
\EP \Big( \Big(X  \frac{\II}{\IIp} \cdot \langle Y,  M \rangle \Big)_{\infty }\Big).
\ede
Furthermore, since $K$ is an $\ff $-predictable process
\bde
 \EP ( (K \cdot ( \Delta Y  \cdot \tJJ ))_{\infty } ) = \EP ( (K \cdot ( \Delta Y \cdot \tJJ )^p)_{\infty } )
 = \EP \Big( \Big(X  \frac{\II}{\IIp} \cdot \check Y \Big)_{\infty }\Big).
\ede
We conclude that for any bounded $\kk$-predictable process $X$
\bde
\EP ( ( X \cdot Y^{\tau })_{\infty }) =
\EP \Big( \Big(X  \frac{\II}{\IIp} \cdot \langle Y,  M \rangle \Big)_{\infty }\Big) +
\EP \Big( \Big(X  \frac{\II}{\IIp} \cdot \check Y \Big)_{\infty }\Big) =
 \EP ( ( X \cdot \wt A )_{\infty }).
\ede
Hence $\wt M = Y^{\tau } - \wt A$ is a $\kk$-local martingale.
It is also clear that $\wt A$ is a $\kk$-predictable process of finite variation.
\end{proof}

\brem
The following question arises: is it possible to derive the equality
$\langle Y , \bar M \rangle =  \langle Y ,  M \rangle + ( \Delta Y \cdot \tJJ)^p$ for any bounded $(\P,\ff)$-martingale $Y$
 using directly the equality $\bar M - M = \tJJ^o - \tJJ^p$?
\erem

For the progressive enlargement $\gg$, we have the following  result.

\bl  \label{lem1.v}
Let $\tau $ be a random time. Then the class of $\gg$-predictable processes is generated by
processes of the form
\bde
X = K \II + g(\tau ) L\JJ = K \I_{\llb 0,\tau\rrb} + g(\tau ) L \I_{\rrb \tau ,\infty \llb}
\ede
where $K$ and $L$ are (elementary) $\ff $-predictable processes and $g: \rr_+ \to \rr$ is a
bounded, Borel measurable function.
If $X$ is bounded then one make take $K$ and $L$ bounded by the same constant.
\el

\ssc{Non-Stopped Processes: Honest Times}

We define the filtration $\gg'$ by the formula
\be \lab{gdash}
\G'_t = \big\{ A \in \G_\infty \,| \,\exists \, \wt A_t, \wh A_t\in \F_t \ \mbox{such that}\
A  = ( \wt A_t \cap\seq{ \tau > t }) \cup ( \wh A_t \cap \seq{\tau \leq t} ) \big\}.
\ee

\bl  \label{lem1.7}
Let $\tau $ be an honest time. Then a process $X$ is a $\gg'$-predictable process if and only if
\bde
X = K \II + L\JJ = K \I_{\llb 0,\tau\rrb} + L \I_{\rrb \tau ,\infty \llb}(t)
\ede
for some $\ff $-predictable processes $K$ and $L$. If $X$ is bounded that one make take $K$ and $L$
bounded by the same constant.
\el

\bl \lab{tt6xx}
Let $X$ be a $\gg'$-predictable process. Then
\be
{}^p (X \II ) = {}^p (K \II )  = {}^p \II K , \quad
{}^p (X \JJ ) = {}^p (L \JJ )  = {}^p \JJ L .
\ee
\el

We are in a position to compute the $(\P,\gg')$-canonical decomposition of a bounded $(\P,\ff)$-martingale $Y$.
We take for granted that $Y$ is a $(\P,\gg')$-semimartingale. The following result can be traced back to
Theorem A in Barlow \cite{B2}, Theorem 5.10 in Jeulin \cite{J2}, Theorem 2 in Jeulin and Yor \cite{JY1}, and Theorem 15 in Jeulin and Yor \cite{JY2}.

\bp \lab{ppp3}
Assume that $Y$ is a bounded $(\P,\ff)$-martingale. The $(\P,\gg')$-canonical decomposition of $Y$
reads $Y  = \widehat M + \widehat A$ where the $\gg'$-predictable process of finite variation $\widehat A$ equals
\be
 \widehat A_t = \Big( \frac{\II}{\IIp} \cdot \langle Y , \bar M \rangle \Big)_t
- \Big( \frac{\JJ}{\JJp} \cdot \langle Y , \bar M \rangle \Big)_t
 =\int_0^t \frac{\II_u}{\IIp_u} \, d \langle Y , \bar M \rangle_u
  - \int_0^t \frac{\JJ_u}{\JJp_u} \, d \langle Y , \bar M \rangle_u .
\ee
\ep

\begin{proof}
Let $X$ be any bounded $\gg'$-predictable process. Then
\bde
\EP ( ( X \cdot Y)_{\infty }) = \EP (( K \II \cdot Y )_{\infty }) + \EP (( L (1-\II ) \cdot Y)_{\infty })
= \EP ( ( K \cdot Y )_{\tau }) - \EP ( ( L \cdot Y )_{\tau })
\ede
since  $\EP (( L \cdot Y)_{\infty })=0$.
By applying  Lemma \ref{rr5} to the bounded $(\P,\ff)$-martingales $K \cdot Y$ and $L \cdot Y$, we obtain
\bde
\EP ( ( K \cdot Y )_{\tau }) = \EP ( \langle K \cdot Y, \bar M \rangle_{\infty }) =
    \EP ( (K \cdot \langle Y, \bar M \rangle )_{\infty })
\ede
and
\bde
\EP ( ( L \cdot Y )_{\tau }) = \EP ( \langle L \cdot Y, \bar M \rangle_{\infty }) =
    \EP ( (L \cdot \langle Y, \bar M \rangle )_{\infty }).
\ede
Lemma \ref{tt6xx} gives
\bde
\EP ( (K \cdot \langle Y, \bar M \rangle )_{\infty }) =
\EP \Big( \Big(K  \frac{\II}{\IIp} \cdot \langle Y, \bar M \rangle \Big)_{\infty }\Big)=
\EP \Big( \Big(X  \frac{\II}{\IIp} \cdot \langle Y, \bar M \rangle \Big)_{\infty }\Big)
\ede
and
\bde
\EP ( (L \cdot \langle Y, \bar M \rangle )_{\infty }) =
\EP \Big( \Big(L  \frac{\JJ }{\JJp} \cdot \langle Y, \bar M \rangle \Big)_{\infty }\Big)=
\EP \Big( \Big(X  \frac{\JJ}{\JJp} \cdot \langle Y, \bar M \rangle \Big)_{\infty }\Big).
\ede
We conclude that for any bounded $\gg'$-predictable process $X$ the following equality holds
\bde
\EP ( ( X \cdot Y)_{\infty }) = \EP ( ( X \cdot \wh A )_{\infty }).
\ede
Hence $\wh M = Y - \wh A$ is a $\gg'$-local martingale.
It is easy to see that $\wh A$ is a $\gg'$-predictable process of finite variation,
so that the equality  $Y  = \widehat M + \widehat A$ is the $(\P,\gg')$-canonical decomposition of $Y$.
\end{proof}

\brem
One can also establish a counterpart of Proposition \ref{ppp2}.
\erem

\end{document}